\documentclass[11pt]{amsart}
\ifx\shlhetal\undefinedcontrolsequence\let\shlhetal\relax\fi

\usepackage{amssymb,amsmath}

\newcommand{\Rightleftarrow}{\Leftrightarrow}
\newcommand{\hra}{\hookrightarrow}
\newcommand{\bbr}{{\mathbb R}}
\newcommand{\bbq}{{\mathbb Q}}
\renewcommand{\lg}{{\rm \ell g\/}}

\newcommand{\bbn}{{\mathbb N}}
\newcommand{\cM}{{\mathcal M}}
\newcommand{\cE}{{\mathcal E}}

\newcommand{\cK}{{\mathcal K}}
\newcommand{\cL}{{\mathcal L}}
\newcommand{\cP}{{\mathcal P}}

\newcommand{\M}{\cM}

\newcommand{\bE}{{\bf E}}
\newcommand{\bc}{{\bf c}}

\newcommand{\w}{{\bf w}}
\newcommand{\bt}{{\mathfrak t}}

\newcommand{\conc}{\char94}
\newcommand{\yes}{{\rm yes}}
\newcommand{\no}{{\rm no}}
\newcommand{\cl}{{\rm cl}}

\newcommand{\GK}{{\mathfrak K}}
\newcommand{\acl}{{\rm acl}}
\newcommand{\prob}{{\rm Prob}}
\newcommand{\imply}{\Rightarrow}
\newcommand{\vare}{\varepsilon}
\newcommand{\Rang}{{\rm Rang}}

\newcommand{\Proof}{\noindent{\sc Proof} \hspace{0.2in}}

\newcommand{\Prob}{{\rm Prob}}

\newcommand{\Dom}{{\rm dom}}

\newcommand{\ex}{{\rm ex}}
\newcommand{\vep}{\varepsilon}

\newcommand{\lesdot}{\mathrel{\mathord{<}\!\!\raise 0.8
pt\hbox{$\scriptstyle\circ$}}}

\newcommand{\rqed}{\hfill\hspace{0.2in}\vrule width 6pt height 6pt depth 0pt
\vspace{0.1in}}

\newtheorem{theorem}{Theorem}[section]
\newtheorem{fact}[theorem]{Fact}

\newtheorem{lemma}[theorem]{Lemma}
\newtheorem{mainlemma}[theorem]{Main Lemma}

\newtheorem{claim}[theorem]{Claim}
\newtheorem{conclusion}[theorem]{Conclusion}
\newtheorem{observation}[theorem]{Observation}

\theoremstyle{definition}

\newtheorem{definition}[theorem]{Definition}

\theoremstyle{remark}

\newtheorem{context}[theorem]{Context}
\newtheorem{remark}[theorem]{Remark}
\newtheorem{general context}[theorem]{General Context}
\newtheorem{notation}[theorem]{Notation}
\newtheorem{discussion}[theorem]{Discussion}

   \def\nofork^#1_#2{\putforkinmargin
	\unionstick^{\textstyle #1}_{\textstyle #2}}

\newbox\noforkbox \newdimen\forklinewidth
\setbox0\hbox{$\textstyle\bigcup$}
\setbox1\hbox to \wd0{\hfil\vrule width \forklinewidth depth \dp0
			height \ht0 \hfil}
\setbox\noforkbox\hbox{\box1\box0\relax}
\def\unionstick{\mathop{\copy\noforkbox}\limits}
\def\nonfork#1#2_#3{#1\unionstick_{\textstyle #3}#2}
\def\nonforkin#1#2_#3^#4{#1\unionstick_{\textstyle #3}^{\textstyle #4}#2}

\title[Zero one laws for graphs... Part I]{Zero one laws for graphs with
edge probabilities decaying with distance. Part I}

\author{Saharon Shelah}
\address{Institute of Mathematics\\
 The Hebrew University of Jerusalem\\
 91904 Jerusalem, Israel\\
 and  Department of Mathematics\\
 Rutgers University\\
 New Brunswick, NJ 08854, USA}
\email{shelah@math.huji.ac.il}
\urladdr{http://www.math.rutgers.edu/$\sim$shelah}
\thanks{The research partially supported by the United States -- Israel
Binational Science Foundation; Publication no 467.}

\begin{document}

\begin{abstract}
Let $G_n$ be the random graph on $[n]=\{1,\ldots,n\}$ with the possible edge
$\{i,j\}$ having probability being $p_{|i-j|}= 1/|i-j|^\alpha$, $\alpha\in
(0,1)$ irrational. We prove that the zero one law (for first order logic)
holds.
\end{abstract}

\maketitle

\section{Introduction}
On 0--1 laws see expository papers e.g., Spencer \cite{Sp}. In {\L}uczak,
Shelah \cite{LuSh:435} the following probabilistic context was
investigated. Let $\bar p=\langle p_i: i\in \bbn\rangle$ be a sequence of
probabilities, i.e. real numbers in the interval $[0,1]_{\bbr}$.  For each
$n$ we draw a graph $G_{n, \bar p}$ with set of nodes $[n]\stackrel{\rm
def}{=}\{1,\ldots, n\}$; for this we make the following independent drawing: 
\begin{itemize}
\item for each (unordered) pair $\{i, j\}$ of numbers from $[n]$ we
draw $\yes$/$\no$ with probabilities $p_{|i-j|}$/ $1-p_{|i-j|}$, and
let
$$
R_n=\{\{i, j\}: i,j \mbox{ are in }[n] \mbox{ and we draw }\yes\}.
$$
\end{itemize}
We consider $R_n$ a symmetric irreflexive 2-place relation. So we have
gotten a random model $\cM^0_{n, \bar p}=([n], R_n)$ (i.e. a graph), but we
also consider the graph expanded by the successor relation $\cM^1_{n, \bar
p}= ([n], S, R_n)$ where $S=\{(\ell, \ell+1): \ell\in\bbn\}$, (more exactly
we use $S_n= S\restriction [n]$), and we may also consider the graph
expanded by the natural order on the natural numbers $\cM^2_{n, \bar p}=
([n], <, R_n)$.  (Here we will give a little background on this structure
below. But the question whether $0-1$ law holds is not discussed here).
Though we shall start dealing generally with random models, the reader can
restrict himself to the case of graphs without losing comprehensibility.

In \cite{LuSh:435} much information was gotten, on when the 0-1 law holds
(see Definition
1.1(1)) and when the convergence law holds (see Definition 1.1(2)), depending
on conditions such as $\sum\limits_{i\in \bbn} p_i<\infty$ and
$\sum\limits_{i\in\bbn} i p_i <\infty$.

The sequences $\bar p$ considered in \cite{LuSh:435} were allowed to be
quite chaotic, and in those circumstances the theorems were shown to be the
best possible, e.g. counterexamples were gotten by replacing $\bar p$ by
$\bar p'$ where for some fast increasing sequence $\langle i_k: k\in
\bbn\rangle$ we let $p'_j= \left\{
\begin{array}{ll} p_k & \ j=i_k \\ 0 & \ (\forall k) j\neq
i_k.\end{array}\right. $

In \cite{Sh:463} a new version of the 0-1 law was introduced, the very weak
zero one law (see 0.1(3), the $h$ variant says that the difference between
the probabilities for $n$ and for $m_n$ when $|n-m_n|\leq h(n)$, converges
to zero) and it was proved for $\cM^2_{n, \bar p}$ when $\sum\limits_i
p_i<\infty$ (we omit $h$ when $h(n)=1$, $m_n=n+1$ and investigate only the
very weak 0-1 law). In \cite{Sh:548} the very weak zero one law was proved
for models with a random two place function and for graphs; Boppana and
Spencer \cite{BoSp} continue this determining the best $h$ for which this
holds. 

Naturally arise the question what occurs if the $p_i$'s are ``well behaved''.
As in Shelah, Spencer \cite{ShSp:304} this leads to considering $p_i=
1/ i^\alpha$ (independently of $n$).
By the results of \cite{LuSh:435}, and (essentially) \cite{ShSp:304}, the
``real'' cases are (on the definition of $\cM^\ell_{n, \bar p}$ see above):
\begin{enumerate}
\item[(A)] $\cM^0_{n, \bar p}$ where $p_i= 1/ i^\alpha$ for $i>1$, $\alpha\in
(0,1)_{\bbr}$ irrational and $p_1=p_2$
\item[(B)] $\cM^1_{n, \bar p}$ where $p_i = 1/ i^\alpha$,
$\alpha\in (0, 1)_{\bbr}$ irrational
\item[(C)] $\cM^2_{n, \bar p}$ where $p_i= 1/ i^\alpha$,
$\alpha\in (1, 2)_{\bbr}$
\end{enumerate}
The main aim of this work is to show that in the case {\bf (A)} we have the
0-1 law, also in case {\bf (B)} we prove the convergence law but at present we
do not know the answer to problem {\bf (C)} (actually analysis indicates
that the problem is whether
there is a formula $\varphi(x)$ which holds in $\cM^2_n$ for $x$ small enough
and fails for $n-x$, $x$ small enough). Here we didn't consider linear
order case. For external reasons the work is divided to two parts, the
second is \cite{Sh:517}.
\noindent {\em Note:} if we let $p_i= 1/ i^{\alpha}$ for $i\geq 1$, surely
$\{\ell, \ell+1\}$ is an
edge, so it is fine, just essentially case {\bf (A)} becomes case {\bf (B)}.
To preserve the distinction between {\bf (A)} and {\bf (B)} we set
$p_1=1/2^\alpha$ in case {\bf (A)}. This is one of many ways to preserve this
distinction; the choice does not matter.

\noindent {\em Main and original context}

\par \noindent
Random graph on $[n]$, with $p_i = 1/i^\alpha$ for $i>1$ and $p_1=p_2$; i.e.
probability of the edge $\{i,j\}$ is $p_{|i-j|}$ and $\alpha\in (0,1)_\bbr
\setminus \bbq$ i.e. is irrational.

But the proofs apply to wider family of cases. We can make a case such that
both \cite{ShSp:304} and \cite{LuSh:435} are particular cases: the probability
for $\{i,j\}$ being an edge of $\cM_n$ for $i,j\in [n]$ is $p^n_{i,j}$. So in
\cite{ShSp:304}, $p^n_{i,j}=p_n$ and in \cite{LuSh:435},
$p^n_{i,j}=p_{|i-j|}$. We can consider $p^n_{i, j}= p^n_{|i-j|}$. We
shall show in another paper
that we shall get the same theory as in case {\bf (A)} above in the limit,
while simplifying the probabilistic arguments, if we change the context to:

\par\noindent
{\em Second context}

\par \noindent
for $\cM_n$ (graph on $\{ 1,\ldots, n\}$) with probability of $\{i,j\}$ being
an edge is $p^n_{i,j}=\frac1{n^\alpha}+\frac1{2^{|i-j|}}$.

So the probability basically has two parts

1) $(\frac{1}{2^{|i-j|}})$:\ \ depends only on the distance, but decays fast,
so the average valency it contributes is bounded.

2) $(\frac{1}{n^\alpha})$:\ \ Does not depend on the distance, locally is
negligible (i.e. for any particular $\{i,j\}$) but has ``large integral''. Its
contribution for the valency of a node $i$ is on the average ``huge'' (still
$\ll n$).

We can think of this as two kinds of edges. The edges of the sort
$n^{-\alpha}$ are as in the paper \cite{ShSp:304}.  The other ones still
give large probability for some $i$ to have valency with no {\em a priori}
bound (though not compared to $n$, e.g.  $\log n$). In this second
context the probability arguments are simpler (getting the same model
theory), but we shall not deal with it here.

\noindent {\em Note:}  If we look at all the intervals $[i, i + k)$, and want
that we get some graph there (i.e. see on H below) and the probability
depends only on $k$ (or at least has a lower bound $>0$ depending only on
$k$), then the chance that for some $i$ we get this graph (by ``second kind
edges'') is $\sim 1$, essentially this behavior stops where $k \approx (\log
(n))^b$ for some appropriate $b > 0$ (there is no real need here to
calculate it).  Now for any graph $H$ on $[k]$ the probability that for a
particular $i<[n-k]$ the mapping $\ell\mapsto i+\ell$ embeds $H$ into
$\cM_n$ is $\geq (\frac{1}{k^\alpha})^{{k}\choose{2}}$ but is $\leq(\frac{1}
{(k/3)^\alpha})^{(k/3)^2}$ (exactly
$$
$$
\noindent $\prod\limits_{\scriptstyle \{\ell,m \} \in {J_1}}
(\frac{1}{|l-m|^\alpha})\cdot\prod\limits_{\scriptstyle \{\ell,m \}
\in {J_2}}(1-\frac{1}{|l-m|^\alpha})\cdot   
p_1^{|\{\ell:(\ell,\ell+1)\ is\ an\ edge \}|}\cdot(1-p_1)^{|\{\ell:
(\ell,\ell+1)\ is\ not\ an\ edge \}|}$
$$
$$
\noindent where $\ell,m \leq k$ and $J_1=\{$ $\{\ell,m \}$: $(\ell,m)$
is an edge and
$|\ell-m|>1 \}$, $J_2=\{$ $\{\ell,m \}$: $(\ell,m)$ is not an edge and
$|\ell-m|>1 \}$. Hence the probability that for no $i<\lfloor n/k\rfloor$
the mapping $\ell\mapsto (k\cdot i +\ell)$ does embed $H$ into $\cM_n$ is
$\leq\bigg(1-\big(\frac{1}{k^{\alpha}}\big)^{{k}\choose{2}}\bigg)^{n/k}$.
Hence if $\beta k^{{\alpha}\cdot{{k}\choose{2}}}=n/k$ that is
$\beta=(\frac{n}{k^{{\alpha}\cdot{{k}\choose{2}}+1}})$
then this probability is $\leq e^{- \beta}$. This is because
$\leq e^{- \beta} \sim \bigg (1- (\frac{\beta}{n})\bigg)^{n}$.
We obtain $(\frac{k \beta}{n}) \leq (\frac{1}{k^{{\alpha}\cdot{{k}
\choose{2}}}})$.
So the probability is small, i.e. $\beta$ large if $k \geq(\frac{2}{\alpha}$
log $n)^{1/2}$; note that the bound for the other direction has the same
order of magnitude. So with parameters, we can interpret, using a sequence
of formulas $\bar\varphi$ and parameter $\bar a$, quite long initial segment
of the arithmetic (see definition below).  This is very unlike
\cite{ShSp:304}, the irrational case, where first order formula
$\varphi(\bar{x})$ really says little on $\bar{x}$: normally it says just
that the $\cl^k$--closure of $\bar x$ is $\bar{x}$ itself or something on
the few elements which are in $\cl^k(\bar{x})$ (so the first order sentences
say not little on the model, but inside a model the first order formula says
little). So this sound more like the $\alpha$ rational case of
\cite{ShSp:304}.  This had seemed like a sure sign of failure of the 0-1
law, but if one goes in this direction one finds it problematic to define
$\bar a_0$ such that $\bar \varphi$ with the parameter $\bar a_0$ defines a
maximal such initial segment of arithmetic, or at least find $\psi(\bar y)$
such that for random enough $\cM_n$, there is $\bar a_0$ satisfying
$\psi(\bar y)$ and if $\bar a_0$ satisfies $\psi(\bar y)$ then $\varphi$
with such parameter define an initial segment of arithmetic of size, say,
$>\log\log\log n$. To interpret an initial segment of arithmetic of size $k$
in $\cM_n$ for $\bar\varphi$ and $\bar a_0$, mean that
$\bar\varphi=\langle\varphi_1 (\bar
x^0,\bar y),\varphi_2 (\bar x^1,\bar y),\varphi_3 (\bar x^2,\bar y)\rangle$ is
a sequence of (first order) formulas, and $\bar a_0$ is a sequence  of
length $\lg(\bar y)$ such that:
the set $\{ x:\cM_n\models\varphi_0(x,\bar a_0) \}$ has $k$ elements, say
$\{ b_0,\ldots,b_{k-1} \}$, satisfying:
$$
\cM_n \models \varphi_1(x_0, x_1, \bar a_0) \Leftrightarrow
\bigvee_{\ell<m<k} (x_0, x_1)= (b_\ell, b_m),
$$
$$
\cM_n \models \varphi_2(x_0, x_1, x_2, \bar a_0) \Leftrightarrow
\bigvee_{\scriptstyle\ell_0,\ell_1,\ell_2>\ell \atop \scriptstyle\ell_2=
\ell_0+\ell_1} (x_0, x_1,x_2)= (b_{\ell_0}, b_{\ell_1}, b_{\ell_2}),
$$
$$
\cM_n \models \varphi_3(x_0, x_1, x_2, \bar a_0) \Leftrightarrow
\bigvee_{\scriptstyle\ell_0,\ell_1,\ell_2< \ell \atop\scriptstyle\ell_2=
\ell_0 \ell_1} (x_0,x_1, x_2) = (b_{\ell_0}, b_{\ell_1}, b_{\ell_2}).
$$
But it is not {\em a priori} clear whether our first order formulas
distinguish between large size and small size in such interpretation.

\noindent{\em Note:} all this does not show why the $0-1$ law holds,
just explain the
situation, and show we cannot prove the theory is too nice (as in
\cite{ShSp:304}) on the one hand but that this is not sufficient for failure
of $0-1$ law on the other hand. Still what we say applies to both contexts,
which shows that results are robust. A nice result would be if we can
characterize $\langle p_i: i\in \bbn \rangle$ such that $\Prob\{i,j\}=p_i
\Rightarrow 0-1$ holds (see below). 

Our idea (to show the $0-1$ law) is that though the ``algebraic closure''
(suitably defined) is not
bounded, it is small and we can show that a first order formula
$\varphi(\bar x)$ is equivalent (in the limit case) to one speaking
on the algebraic closure of $\bar x$.

Model theoretically we do not get in the limit a first order theory which is
stable and generally ``low in the stability hierarchy'', see Baldwin, Shelah
\cite{BlSh:528}, for cases with probability $\sim n^{-\alpha}$ (the
reason is of course that restricted to ``small'' formulas in some
cases there is a definable linear order (or worse)). However we
get a variant of stability over a predicate: on ``small'' definable sets the
theory is complicated, but for types with no small formulas we are in the
stable situation. In fact the model theoretic setting is similar to the one
in \cite{Sh:463} but we shall not pursue this.

Note that Baldwin, Shelah \cite{BlSh:528} deal with random models with more
relations $R$ with probabilities $n^{\alpha(R)}$ (satisfying the parallel to
irrationality of $\alpha$). There, the almost sure theory is stable. In
\cite{Sh:550} we define a family of 0-1 contexts where further drawings of
relations give us a new context in this family and in all such contexts,
elimination of quantifiers to the algebraic closure (as in \cite{ShSp:304},
\cite{BlSh:528}) holds, but the context is possibly ``almost nice''
not nice, i.e. we allow that
every $\bar a$ has a nontrivial closure, as in the case in which we have the
successor function. Here this is dealt within the general
treatment of the elimination, but not used in the main case $\cM^0_n$. We
could deal with abstract version allowing further drawing as in \cite
{Sh:550} also here.

See more \cite{Bl96}, \cite{Sh:637}.

We have chosen here quite extreme interpretation of ``$\bar p$ is simple,
simply defined''. It seems desirable to investigate the problem under more
lenient conditions. A natural such family of $\bar p$'s is the family of
monotonic ones. Can we in this family characterize
$$
\{\bar p:\bar p\mbox{ monotone sequence, }\cM^0_{n,\bar p}\mbox{ satisfies
the }0-1
\mbox{ law}\}?
$$
This will be addressed and solved in \cite{Sh:581}.

The two cases considered above are protypes of some families with the 0-1 law,
but there are some others, for example with the value of the exponent $\alpha$
``in the appropriate neighbourhood'' of a rational (and some degenerate
ones of course).

Let us review the paper.

\noindent{\em Note:} in \S1 -- \S3 we deal with general contexts. In these
three sections sufficient conditions are proven for the 0-1 law to hold in
0-1 context; for notational simplicity we restrict ourselves to vocabulary
which
contains finitely many predicate relations (not only a symmetric irreflexive
2-place relation). The proof is based on elimination of quantifiers by the
help of the closure without using probability arguments. Note that in the
application we have in mind, the closure has order of magnitude up to
$\sim$ log$|\cM_n|$. In \cite{ShSp:304} $\cl$ is bounded i.e. $|\cl(A)|$ has a
bound depending on $|A|$ (and $\alpha$ of course) only while here is not
bounded. In the second part, \S4
-- \S6 deal with $\cM^0_n$ and \S7 deal with $\cM^1_n$.

In \S1 we give the basic definitions, including $A<_i
B$ (intended to mean: $B$ is the algebraic closure of $A$ but this closure has
no {\em a priori} bound). The restriction to: $\cM_n$ has set of elements $[n]$
(rather than some finite set) is not important for the proof. In \S1,
$A<_i B$
and $A<_s B$ are defined in terms of the number of embeddings of $A$ into
$\cM_n$ in a sufficiently random model, and from $<_i$ we define
$\cl^k(A,M)$.

In \S2 a fundamental relation (i.e. given a priori) on structures $M$ is
$\cl^k$. From it notions of
$A<_i B$ and $A<_s B$ are defined in terms of embeddings $f\subseteq g$ of
$A,B$ into a sufficiently random $\cM_n$ and the relations between
$g(B)$ and $\cl^k_{\cM_n}(f(A), \cM_n)$. Then these definitions are
reconciled with those in \S1, when the closure is chosen as in \S1.
Two axiomatic frameworks for an abstract elimination of quantifiers argument
are presented. (This generalizes \cite{BlSh:528}.) These frameworks and
further conditions on $\cl^k$ provide sufficient conditions for $0-1$ laws
and convergence laws. 

\noindent{\em Note:} in \S2 we retain using ``relation free amalgamation'' (as
in \cite{BlSh:528}, but in \cite{Sh:550} we will use more general one).
However we waive ``random $A$ has no non-trivial closure'', hence use ``almost
nice'' rather than nice (and also waive the {\em a priori} bounds on closure).

In \S3 we deal with the case where the natural elimination of quantifiers is to
monadic logic. This seems natural, although it is not used later.

We now proceed to describe part II, the main point of \S4 is to introduce a
notion of weight ${\bf w}(A,B,\lambda)$ which depends on an equivalence
relation $\lambda$ on $B\setminus A$.
(Eventually such $\lambda$ will be defined in terms of the ``closeness'' of
images of points in $B$ under embeddings into $\cM_n$.) Relations $A\leq_i^*
B$ and $A\leq_s^* B$ are defined in terms of ${\bf w}$. The intension is that
$\leq_i^*$ is $\leq_i$ etc, thus we will have direct characterization of the
later.

\S5 contains the major probability estimates. The appropriate $\lambda$ is
defined and thus the interpretations of $<_i^*$ and $<_s^*$ in the first
context $(\cM^0_n, p_i=\frac{1}{i^\alpha})$. Several proofs are analogous to
those in \cite{ShSp:304} and \cite{BlSh:528}, so we treat them only
briefly. The new point is the dependence on distance, and hence the
equivalence relations $\lambda$.

In \S6 it is shown that the $<^*_i$ and $<^*_s$ of \S5 agree with the $<_i$
and $<_s$ of \S1. Further, if $\cl^k$ is defined from the weight function in
\S4, these agree with $<_i$, $<_s$ as in \S2 and we prove the ``simple almost
niceness'' of Definition \ref{2.5B}, so the 'elimination of quantifiers
modulo quantification on (our) algebraic closure' result
applies. This completes the proof of the $0-1$ law for the first
context. The model theoretic considerations in the proof of this
version of niceness (e.g. the compactness) were less easy than I expect.

\S7 deals with the changes needed for $\cM^1_{n,\bar{p}}$ where only the
convergence law is proved.

\noindent
{\em Note:} our choice ``$\cM_n$ has set of element $[n]$'' is just for
simplicity (and tradition), we could have $\cM_n$ has set of elements a finite
set (not even fixed) and replace $n^\varepsilon$ by $\|\cM_n\|^\varepsilon$ as
long as ``for each $k$ for every random enough $\cM_n$ we have
$\|\cM_n\|>k$''. Also the choice of $n^\varepsilon$ in Definition 1.2 is the
most natural but not unique case. The paper is essentially self contained,
assuming only basic knowledge of first order logics and probability.

\begin{notation}
\begin{itemize}
\item $\bbn$ is the set of natural numbers ($\{0, 1,2,\ldots\}$)
\item $\bbr$ is the set of reals
\item $\bbq$ is the set of rationals
\item $i, j, k, \ell, m, n, r, s, t$ are natural numbers and
\item $p$, $q$ are probabilities
\item $\alpha$, $\beta$, $\gamma$, $\delta$ are reals
\item $\varepsilon$, $\zeta$, $\xi$ are positive reals (usually quite small)
and also $\bc$ (for constant in inequalities)
\item $\lambda$ is an equivalence relation
\item $M, N, A, B, C, D$ are graphs or more generally models (
that is structures, finite of fixed finite vocabulary, for notational
simplicity with predicates only, if not said otherwise; the reader can
restrict himself to graphs) 
\item $|M|$ is the set of nodes or elements of $M$, so $\|M\|$ is the number
of elements.
\item ${ \mathcal M }$ denotes a random model,
\item  $ \mu $ denotes a distribution (in the probability sense),
\item $[n]$ is $\{1, \ldots, n\}$
\item $A\subseteq B$ means $A$ is a submodel of $B$ i.e. $A$ is $B$ restricted
to the set of elements of $A$ (for graphs: induced subgraph)\\
We shall not always distinguish strictly between a model and its set of
elements. If $X$ is a set of elements of $M$, $M\restriction X$ is $M$
restricted to $X$.
\item $a$, $b$, $c$, $d$ are nodes of graphs / elements of models
\item $\bar a$, $\bar b$, $\bar c$, $\bar d$ are finite sequences of
nodes / elements
\item $x$, $y$, $z$ variables
\item $\bar x$, $\bar y$, $\bar z$ are finite sequences of variables
\item $X$, $Y$, $Z$ are sets of elements
\item $\tau$ is a vocabulary for simplicity with predicates only (we may
restrict a predicate to being symmetric and/or irreflexive (as for
graphs)),
\item $\cK$ is a family of models of fixed vocabulary, usually
$\tau=\tau_\cK$
\item the vocabulary of a model M is  $\tau_M$,
\item $\bar a \conc \bar b$ or $\bar a \bar b$ is the concatenation of
the two sequences, $\bar a\conc b$ or $\bar a b$ is $\bar
a\conc\langle b\rangle$
\item the extensions $g_1$, $g_2$ of $f$ are disjoint if $x_1\in
\Dom(g_\ell)\setminus \Dom(f)$, $x_2\in \Dom(g_{3-\ell}) \Rightarrow
x_1\neq x_2$.
\end{itemize}
\end{notation}
\medskip

\noindent{\bf Acknowledgements:}\quad We thank John Baldwin and Shmuel
Lifsches and \c{C}i\u{g}dem Gencer and Alon Siton for helping in various
ways and stages to make the paper more user friendly.

\section{Weakly nice classes}
We interpret here ``few'' by: ``for each $\varepsilon$ for every random enough
$\cM_n$, there are (for each parameter) $<n^\varepsilon$''. We could use other
functions as well.

\begin{general context}
\label{1.0}

(i) Let $\tau$ be fixed vocabulary which for simplicity having only
predicates, i.e. symbols for relations.

(ii) $\cK$ be a class of finite $\tau$- models closed under isomorphism
and submodels. For $n \in \bbn$, $\cK_n$ is a set of $\tau$-models
which usually have universe $[n]=\{1,...,n \}$(just for notational
simplicity).

(iii) Let ${\cM}_n$ be a random model in a fixed vocabulary $\tau$
which is an element of $\cK_n$, that is we have $\mu_n$ a function
such that $\mu_n:\cK_n \rightarrow [0,1]_{\bbr}$ and
$\sum \{ \mu_n (\cM):\cM \in \cK_n \}=1$, so $\mu_n$ is called
a distribution and $\cM_n$
the random model for $\mu_n$, so we restrict ourselves to finite or
countable $\cK_n$. We omit $\mu_n$ when clear from the context.

(iv) We call $(\cK,\langle(\cK_n,\mu_n): n<\omega\rangle)$  a $0-1$ context
and denote it by $\GK$ and usually consider it fixed; we may 'forget' to
mention $\cK$. So ,

(v) The probability of $\cM_n\models \varphi$; $Prob(\cM_n\models \varphi)$ is
\[\sum \{\mu_n(\cM): \cM  \in \cK_n, \cM \models \varphi \}.\]

(vi) The meaning of ``for every random enough $\cM_n$ we have $\Psi$'' is
$$
\langle\Prob (\cM_n \models \Psi): n<\omega\rangle\ \mbox{ converges to 1;}
$$
alternatively, we may write ``almost surely $\cM_n\models \Psi$''.

(vii) We call $\GK$ a $0-1$ context if it is as above.
\end{general context}

\begin{definition}
\label{1.1}
\begin{enumerate}
\item The $0-1$ law (for $\GK$) says: whenever $\varphi$ is a f.o.
(=first order) sentence in vocabulary $\tau$,
$$
\langle \Prob((\cK_n,\mu_n): n<\omega\rangle) n<\omega\rangle\mbox{
converges to } 0
\mbox{ or \ to } 1.
$$
\item The convergence law says: whenever $\varphi$ is a f.o. sentence in
$\tau$,
$$
\langle \Prob(\cM_n\models \varphi): n<\omega\rangle\mbox{ is a convergent
sequence}.
$$
\item The very weak $0-1$ law says: whenever $\varphi$ is a f.o. sentence in
$\tau$,
$$
\lim\limits_n[\Prob(\cM_{n+1}\models \varphi)-\Prob(\cM_n\models\varphi)]=0.
$$
\item The $h$-very weak $0-1$ law for $h:\bbn\rightarrow\bbn\setminus \{0\}$
say: whenever $\varphi$ is a f.o sentence in $\tau$,
$$
0 =\lim_n\max_{\ell, k\in [0, h(n)]}|
\Prob(\cM_{n+k}\models\varphi)-\Prob(\cM_{n+
\ell}\models\varphi)|
$$
\end{enumerate}
\end{definition}

\begin{notation}
\label{1.1A}
$f:A\hra B$ means: $f$ is an embedding of $A$ into $B$ (in the model theoretic
sense, for graphs: isomorphism onto the induced subgraph).
\end{notation}

\begin{definition}
\label{1.2}
\begin{enumerate}
\item Let
$$
\begin{array}{ll}
{\cK}_\infty =\Big\{A:&A\mbox{ is a finite $\tau$-model}\\
\ &0<\lim\limits_n\sup[\Prob((\exists f)
(f: A\hra \cM_n))]\Big\}
\end{array}
$$
recall (1.1(v)) that Prob(($\exists f)(f: A\hra \cM_n))=
\sum \{\mu_n(\cM_n): \cM \in \cK_n$ and there is an embedding
$f: A\hra \cM_n \}$, $n< \omega$.

Also let $T_{\infty} =^{\rm df }\{\varphi:\varphi$ is a f.o. sentence
in the vocabulary of $ {\mathcal K} $ such that every
random enough $ { \mathcal M } _n $ satisfies it $\}$.

\item $A\leq B$ means: $A,B \in {\cK}_{\infty}$ and $A$ is
a submodel of $B$; of course $A<B$ means $A\leq B$ and
$A\neq B$, similarly for others below.
\item $A\leq_i B$ means: $A\leq B$ and for each
$\varepsilon \in \bbr^+$ we
have:
$$
1=\lim\limits_n\left[\Prob\left(
\begin{array}{l}
\mbox{if } f_0:A\hra \cM_n\\
\mbox{then the number of $f_1$ satisfying }\\
f_0\subseteq f_1:B\hra \cM_n\mbox{ is }\leq n^\varepsilon.
\end{array}
\right)\right]
$$
Also let $\ex( f_0 , B, M )=\ex (f_0, A,B,M)=^{\rm df}\{f:f$ is an embedding
of $B$ into $M$ extending $f_0\}$. 
\item $A \leq_s B$ means: $A \leq B$ and there is no $C$ with $A<_i C \leq B$.
\item $A <_{pr} B$ means: $A <_s B$ and there is no $C$ with $A<_s C<_s B$
($pr$ abbreviates {\em primitive}).
\item $A <_a B$ means that $A\leq B$ and, for some $\varepsilon\in\bbr^+$ for
every random enough $\cM_n$, for no $f: A\hra \cM_n$ do we have
$n^\varepsilon$ pairwise disjoint extensions $g$ of $f$ satisfying $g: B\hra
\cM_n$.
\item $A\leq_{m}^{s} B$ means $A\subseteq B$ are from $\cK$ and for every
$X\subseteq B$ with $\leq m$ elements, we have $A\restriction(A\cap X)\leq_s
(B\restriction X)$.
\item $A\leq_{k, m}^{i} B$ means $A\subseteq B$ are from $\cK$ and for every
$X\subseteq B$ with $\leq k$ elements there is $Y$, $X\subseteq Y\subseteq B$
with $\leq m$ elements such that $A\restriction (A\cap Y)\leq_i (B\restriction
Y)$.
\item For $h: \bbn\times \bbr^+ \rightarrow \bbr^+$, we define $A\leq^h_i B$
as in part (3) replacing $n^\varepsilon$ by $h(n,\varepsilon)$, and similarly
$A\leq^h_a B$ (in part (6)), hence $A\leq^h_s B$,
$A\leq^h_{pr} B$, $A<^h_a B$,
$A\leq^{s,h}_m B$, $A\leq^{i, h}_{k, m} B$.
\end{enumerate}
\end{definition}

\begin{remark}
\label{1.2A}
\begin{enumerate}
\item In these circumstances the original notion of algebraic closure
is not well behaved. $A \leq_i B$ provides a reasonable substitute for
$A\subseteq B\subseteq \acl (A)$.\\
\item Note: for $\leq^h_i$ to be transitive we need: for every
$\varepsilon_1>0$ for some $\varepsilon_2>0$ for every $n$ large
enough $h(n, \varepsilon_2)\times h(n, \varepsilon_2)\leq h(n,
\varepsilon_1)$.
\item Why do we restrict ourselves to $\cK_\infty$ (in \ref{1.2}(1)-(6))? The
relations in \ref{1.2}(1)-(6) describe situation in the limit. So why in
\ref{1.2}(7), (8) do we not restrict ourselves to $A,B\in\cK_\infty$? As for
$A\in \cK_\infty$, for quite random $\cM_n$, and $f:A\hra\cM_n$ the set
$\cl^k(f(A), \cM_n)$ may be quite large, say with $\log(n)$ elements, so it
(more exactly the restriction of $\cM_n$ to it) is not necessarily in
$\cK_\infty$; this is a major point here.
\end{enumerate}
\end{remark}

Let us expand.\\
If $A\in\cK$ has a copy in a random enough $\cM_n$ and we have $0-1$ law then
$T_\infty$(see 1.4(1)) says that copies of $A$ occur. But if $\cM_n$ is
random enough, and for example
$A=\{a_1,a_2,a_3\}\leq\cM_n$, and $B=\cM_n\restriction \cl^k(\{a_1,a_2,a_3\},
\cM_n)$ has, say, $\log(n)$ elements {\em then} it does not follow that
$T_\infty\models$``a copy of $B$ occurs'', as $\cM_n$ may not be random enough
for $B$. Still for the statements like
$$
(\exists x_1,x_2,x_3)(\cl^k(\{x_1,x_2,x_3\})\models \varphi)
$$
the model $\cM_n$ may be random enough. The point is that the size of $B$
could be computed only after we have $\cM_n$.\\
Another way to look at it: models $M_\infty$ of $T_\infty$ are very random in
a sense, but $\cl(\{a_1,a_2,a_3\},M_\infty)$ is infinite, may even be
uncountable, so
randomness concerning it becomes meaningless.

\begin{definition}
\label{1.3}
For $A\subseteq M$ and $k<\omega$ define
\begin{enumerate}
\item[(a)] $\cl^k(A, M)=\bigcup\{B: B\subseteq M,\; B\cap A\leq_i B,\mbox{ and
}|B|\leq k\}$,
\item[(b)] $\cl^{k,0}(A,M)=A$,
\item[(c)] $\cl^{k,m+1}(A,M)=\cl^k(\cl^{k,m}(A,M),M)$.
\end{enumerate}
\end{definition}

\begin{observation}
\label{1.4}
1) For all $\ell,k\in \bbn$ and $\vare\in \bbr^+$ we have
$$
1=\lim\limits_n\big[\prob\big(A\subseteq\cM_n,|A|\le\ell\imply
|\cl^k(A,\cM_n)|<n^\vare\big)\big].
$$
2) Moreover, for every $k\in\bbn$ and $\varepsilon\in\bbr^+$ for some
$\zeta\in \bbr^+$ (actually, any $\zeta<\varepsilon/(k+1)$ will do) we have
$$
1=\lim_n\big[\prob(|A|\leq\cM_n, |A|\leq n^\zeta\Rightarrow |\cl^k(A,\cM_n)| <
n^\varepsilon)\big].
$$
\end{observation}

\begin{remark}
True for $\cl^{k, m}$ too, but we can use claim \ref{1.11}
instead.
\end{remark}

\begin{definition}
\label{1.5}
$\GK=\langle\cM_n:n<\omega\rangle$ is {\em weakly nice} if whenever $A<_s B$
(so $A\neq B$), there is $\varepsilon\in\bbr^+$ with
$$
1=\lim\limits_n\left[\prob\left(
\begin{array}{l}
\mbox{if } f_0:A \hra \cM_n\ \mbox{ then there is } F \mbox{ with }
|F|\ge n^\varepsilon\ \mbox{ and }\\
\mbox{ (i) } f\in F\imply f_0 \subseteq f: B \hra \cM_n\\
\mbox{ (ii) } f'\neq f''\in F\imply\Rang(f')\cap\Rang(f'')=\Rang (f_0)
\end{array}
\right)\right].
$$
If clause (ii) holds we say the $f\in F$ are pairwise disjoint over $f_0$ or
over $A$. In such circumstances we say that $\varepsilon$ witnesses $A<_s B$.
\end{definition}

\begin{remark}
Being weakly nice means there is a gap between being pseudo algebraic
and non-pseudo algebraic (both in our sense), so we have a strong dichotomy.
\end{remark}

\begin{fact}
\label{1.6}
For every $A, B, C$ in ${\cK} {}_\infty$:
\begin{enumerate}
\item $A \le_i A$,
\item $A \le_i B, \; B \le_i C \Rightarrow A \le_i C$,
\item $A \le_s A$,
\item if $A_1\leq B_1$, $A_2\leq B_2$, $A_1\leq A_2$, $B_1\leq B_2$, $B_1
\setminus A_1=B_2\setminus A_2$ then  $A_2\leq_s B_2\Rightarrow A_1\leq_s B_1$
and $A_1\leq_i B_1 \Rightarrow A_2 \leq_i B_2$,
\item $A<_i B$ iff for every $C$ we have $A\leq C<B\ \Rightarrow\ C<_a B$.
\end{enumerate}
\end{fact}

\Proof Easy (e.g. \ref{1.6}(5) by the $\Delta$-system argument
(for fixed size of the sets and many of them); note $|B|$ is
constant).
\hfill\rqed$_{\ref{1.6}}$

\begin{claim}
\label{1.7}
If $A <_s B <_s C$ then $A <_s C$
\end{claim}

\Proof
{\em First proof:}

If not, then for some $B'$ we have $A<_i B'\leq C$. If $B'\subseteq B$
we get contradiction to $A<_s B$, so assume $B'\nsubseteq B$. By
\ref{1.8}(1) below we have $(B'\cap B)<_i B'$ so by \ref{1.6}(4) we
have $B<_i (B\cup B')$, hence we get contradiction to $B<_s C$.

\noindent
{\em Second proof:}
(Assuming $\GK$ is weakly nice i.e. if we define $<_s$ by \ref{1.5}.)
Let $\vare>0$ witness $A<_s B$ in Definition \ref{1.5} and let
$\zeta>0$ witness $B<_s C$ in Definition \ref{1.5}. Choose $\xi=
\min\{\varepsilon/2,\zeta/2\}$; (actually just $\xi<\vep \wedge \xi<
\zeta$ suffice). Let $n$ be large enough; in particular
$n^\varepsilon>|C|$ and let $f_0: A\hra\cM_n$. So we have (almost surely)
$\{f^i_1: i<i^\ast\}$, where $i^\ast\ge n^\varepsilon$, and $f_0\subseteq
f^i_1$ and $f^i_1: B\hra\cM_n$ and the $f^i_1$'s are pairwise disjoint
over $A$.

Now, almost surely for every $i$ we have $\{f^{i,j}_2: j<j^*_i\}$ with $f^i_1
\subseteq f^{i,j}_2$ and $f^{i,j}_2:C\hra\cM_n$ and, fixing $i$, the
$f^{i,j}_2$'s are pairwise disjoint over $B$ and $j^*_i\geq n^\zeta$.\\
Clearly (when the above holds) for $\ell^\ast=n^\xi$ we can find $\{j_k:
k\leq\ell^\ast\}$ such that $\{f_2^{k,j_k}:k<\ell^\ast\}$ are pairwise
disjoint over $A$ (just choose $j_k$ by induction on $k$ such that:
$\Rang(f_2^{k,j_k}\restriction(C\setminus B))$ is disjoint to
$$
\bigcup\{\Rang(f^{i}_1\restriction (B\setminus A)):i< \ell^*\}\cup\bigcup
\{\Rang(f^{i,j_i}_2\restriction (C\setminus B)): i<k\};
$$
at stage $k$, the number of inappropriate $j<n^\zeta$ is
$$
\leq |C\setminus B|\times k+|B\setminus A|\times\ell^*\leq |C|\times\ell^*=
|C|\times n^\xi).
$$
\hfill\rqed$_{\ref{1.7}}$

\begin{fact}
\label{1.8}
Suppose $A \le B \le C$.
\begin{enumerate}
\item If $A\le_i C$ {\em then} $B \le_i C$.
\item If $A\le_s C$ {\em then} $A\le_s B$.
\item If $A<_{pr}C$ and $A\leq_s B\le_s C$ {\em then} either $B=A$ or $B=C$.
\end{enumerate}
\end{fact}

\Proof Reread the definitions.

\begin{fact}
\label{1.9}
\begin{enumerate}
\item If $A\le_s B$ {\em then} there is some $n<\omega$ and a sequence
$\langle A_l: l\leq n \rangle$ such that $A=A_0<_{pr} A_1 <_{pr}\ldots<_{pr}
A_n= B$ (possibly $n=0$).
\item If $A <_{pr} C$ and $A < B < C$ then $B <_i C$
\end{enumerate}
\end{fact}

\Proof  For proving (2), choose a maximal $B'$ such that $B\le_i B'\le C$, it
exists as $C$ is finite (being in $\cK_\infty$), and as $B\leq_i B$ (by
\ref{1.6}(1)). It follows that if $B'< B^{''}\le C$ then $\neg B'\le_i B^{''}$
(by \ref{1.6}(2)). Hence $B'\le_s C$. But $A<_{pr} C$ hence by the Definition
\ref{1.2}(5) we have $A<_s C$ so by \ref{1.8}(2) $A<_s B'$; so by the
definition of $<_{pr}$ we have $B'=C$, so $B\leq_i B'=C$ as required. Part (1)
is clear as $C$ is finite (being in $\cK_\infty$) and the definition of
$\leq_{pr}$. \hfill\rqed$_{\ref{1.9}}$

\begin{claim}
\label{1.10}
${\cK}$ is weakly nice iff whenever $A<_{pr} B$ there is $\vare\in\bbr^+$ such
that
$$
1=\lim\limits_n\left[\prob\left(
\begin{array}{l}
\mbox{if }\ f_0: A\hra\cM_n\mbox{ then there is } F\mbox{ with } |F|\ge
n^\varepsilon \mbox{ and}\\
\qquad f_1 \in F \ \Rightarrow\  f_0 \subseteq f_1 : B\hra\cM_n
\end{array}
\right)\right]
$$
\end{claim}

\Proof  $\Rightarrow$ is obvious (as $A <_{pr} B$ implies $A <_s B$).

Let us prove $\Leftarrow$: we have $A\le_s B$ and by fact \ref{1.9}(1) there
is a sequence $A=A_0<_{pr} A_1<_{pr}\cdots<_{pr} A_k=B$. The proof is by
induction on $k$. The induction step for $k>1$ is by the second proof of
\ref{1.7} and $k=0$ is
\ref{1.6}(3). So assume $k=1$, hence $A<_{pr} B$. By fact \ref{1.9}(2) if
$A<B'\le B$ then $B'\le_i B$. Fix $p\in (0,1)_{\bbr}$. If $n$ is large enough
then the probability of having both
\begin{enumerate}
\item[(a)] for every $f_0:A\hra\cM_n$ there are at least $n^\varepsilon$
different extensions $f^i_1$ satisfying $f_0 \subseteq f^i_1 : B\hra\cM_n$ and
\item[(b)] for every $a\in B\setminus A$ and $f^+_0: A\cup \{a\}\hra \cM_n$
there are at most $n^{\varepsilon/2}$ different extensions $f^i_2$ satisfying
$f^+_0 \subseteq f^i_2: B\hra\cM_n$
\end{enumerate}
is $\ge 1-p$ (for clause (b) use $A\cup\{a\}<_i B$ for every $a\in B\setminus
A$ which holds by \ref{1.9}(2)). Let $f_0:A\hra\cM_n$, and let $\langle f^j_1:
j<j^\ast\rangle$ be a maximal family of pairwise disjoint extensions of $f_0$
to an embedding of $B$ into $\cM_n$. Let $F=\{f:$ $f$ is an embedding of $B$
 into
$\cM_n$ extending $f_0 \}$. By (b) we have
$$
n^\varepsilon\leq |F| \le j^\ast\times|B\setminus A|\times|B\setminus A|\times
n^{\varepsilon/2}.
$$
Hence if $n$ is large enough, $j^\ast>n^{\varepsilon/3}$ (with probability
$\ge 1-p$), and this is enough. \hfill\rqed$_{\ref{1.10}}$

\begin{claim}
\label{1.11}
$\cl^{k,m}(A,M)\subseteq\cl^{k^\ast}(A, M)$ where $k^\ast= k^m$.
\end{claim}

\Proof Define $k(\ell)$ by induction on $\ell\leq m$: $k(0)=1$, $k(1)=k$
and for $\ell<m$ (but $\ell\geq 1$),
$k(\ell+1):=k(\ell)k$. For $\ell\le m$ define
$A_\ell=\cl^{k,\ell}(A,M)$. Now if $x\in A_m$ then there is some $\ell<m$ such
that $x\in A_{\ell+1}\setminus A_\ell$.  Let us prove by induction on $\ell\le
m$ that $x\in A_\ell\Rightarrow x\in\cl^{k(\ell)}(A,M)$. For $\ell=0$
and $\ell=1$ this is
clear. If $x\in A_{\ell+1}\setminus A_\ell$ then there is $C$ with $|C|\le k$
such that $x\in C$ and $C\cap A_\ell<_i C$. By the induction hypothesis, for
$y\in C\cap A_\ell$ we have $y\in\cl^{k(\ell)}(A, M)$ hence there is $C_y$
with $|C_y|\le k(\ell)$ such that $y\in C_y$ and $C_y\cap A<_i C_y$. Let $C^0
=\bigcup\limits_{y\in C\cap A_\ell} C_y \cap A$, $C^1=\bigcup\limits_{y\in
C \cap A_\ell} C_y$ and $C^2=C^1\cup C$. As $|C|\leq k$, we get
$$
|C^2|\le k(\ell)\cdot|C\cap A_\ell|+|C\setminus A_\ell|\le k(\ell)\cdot k
\le k(\ell+1),
$$
so (as $x\in C^2$) it suffices to show that $C^0\le_i C^2$ and by transitivity
(i.e. by \ref{1.6}(2)) it suffices to show that $C^0\le_i C^1$ and that $C^1
\le_i C^2$. Why $C^1\le_i C^2$? Because $C\cap A_\ell\le_i C$ and $C\cap
A_\ell \subseteq C^1\subseteq A_\ell$ and hence $C^1\le_i C^1\cup C= C^2$ by
\ref{1.6}(4). Why $C^0\le_i C^1$? Let $C\cap A_\ell=\{y_s:s<r\}$. Now $C^0
\le_i C^0\cup C_{y_0}$ by \ref{1.6}(4) because $A\cap C_{y_0}\le_i
C_{y_0}$ and $A\cap C_{y_0}\subseteq C^0$ and similarly by induction
$$
C^0 \le_i C^0 \cup C_{y_0}\le_i C^0\cup C_{y_0}\cup C_{y_1}\le_i\ldots\le_i
C^0\cup\bigcup\limits_{s<r} C_{y_s} = C^1.
$$
So as $\leq_i$ is transitive (\ref{1.6}(2)) we are done.
\hfill\rqed$_{\ref{1.11}}$

\begin{claim}
\label{1.12}
For every $\varepsilon\in\bbr^+$ and $\ell,k,m$ we have
$$
1=\lim\limits_n\left[\prob\left(
\begin{array}{l}
\mbox{ if }A\in\cK_\infty,|A|\le\ell\mbox{ and }f:A\hra\cM_n\\
\mbox{ then }|\cl^{k,m}(f(A),\cM_n)|<n^\varepsilon\\
\end{array}
\right)\right].
$$
\end{claim}

\Proof  By the previous claim \ref{1.11}, w.l.o.g. $m =1$. This holds by
Definition \ref{1.2}(3) and Definition \ref{1.3}. \hfill\rqed$_{\ref{1.12}}$

\begin{fact}
\label{1.13}
\begin{enumerate}
\item For every $A$ and $m, k$, for any $M\in \cK$ if $f:A\hra M$ then
\begin{enumerate}
\item[$(\alpha)$] $\cl^{k,m}(f(A),M)\leq^i_{1,k}\cl^{k,m+1}(f(A),M)$,
\item[$(\beta)$] for some $m'=m'(k, m)$ we have
$$
f(A)\leq^i_{k, m'}\cl^{k,m}(f(A),M)
$$
(we can get more),
\item[$(\gamma)$] $f(A)\le_i\cl^{k,m}(f(A),M_n)$ or the second is not in
$\cK_\infty$.
\end{enumerate}
\item For every $m$, $k$, $\ell$ for some $r$ we have:

for any $A\in\cK_\infty$,
$$
1=\lim_n\left[\Prob\left(\mbox{if } f:A\hra\cM_n\mbox{ then }
f(A)\leq^i_{\ell,r}\cl^{k,m}(f(A),\cM_n)\right)\right].
$$
\end{enumerate}
\end{fact}

\begin{remark}
In our main case $\cK=\cK_\infty$.\\
Recall for \ref{1.13}(1)$(\gamma)$ that $\cl^{k,m}(f(A),\cM_n)$ is in
general not necessarily in $\cK_\infty$.
\end{remark}

\Proof  1) We leave the proof of $(\alpha)$ and $(\beta)$ to the reader.
For proving clause ($\gamma$),
let $A_0=f(A)$ and for $\ell\le n$ let $A_\ell=
\cl^{k,\ell}(f(A),M)$, and assume $A_n\in\cK_\infty$. So for $\ell<m$ we have
$A_{\ell+1}=A_\ell \cup\bigcup\limits_{j<m_\ell} C_{\ell, j}$ with
$|C_{\ell,j}| \le k$ and $A_{\ell+1}\cap C_{\ell,j}\leq_i C_{\ell,j}$.
It follows by
\ref{1.6}(4) that $\langle A_\ell\cup\bigcup\limits_{i<j} C_{\ell,i}: j\le
m_\ell\rangle$ is $\leq_i$-increasing and $A_\ell\leq_i A_{\ell+1}$. By
induction we get $A_0 \leq_i A_m$ which is the desired conclusion.

\noindent 2) Read the proofs of \ref{1.13}(1) + \ref{1.11}.
\hfill\rqed$_{\ref{1.13}}$

\begin{remark}
In a more general context the previous conclusion is part of the definition of
``${\cK}$ is nice'' and also $\nonforkin{}{}_{}^{}$ of \ref{1.16}
below is a basic property (on the later see \cite{Sh:550}).
\end{remark}

\begin{fact}
\label{1.14}
${\cK}_\infty$ is closed under isomorphisms and taking submodels.
\end{fact}

\begin{fact}
\label{1.15}
For every $\ell, k, m$ there is a first order formula
$\varphi(y,x_0,\ldots,x_{\ell-1})$ such that for every $M\in{\cK}$ and
$b,a_0,\ldots,a_{\ell-1}$ in $M$ for any $f: A \hra \cM$
$$
M\models\varphi(b,a_0,\cdots,a_{\ell-1})\quad\mbox{ iff }\quad b\in\cl^{k,m} (
\{a_0,\cdots,a_{\ell-1}\},M).
$$
\end{fact}

\Proof  By finiteness of $\tau$ (as $\tau_{\cK}$ is having no function
symbols); or see proof of clause $(\beta)$ of $2.6$.

\begin{definition}
\label{1.16}
$\nonforkin{C_1}{C_2}_{B}^{D}$ means: they are all submodels of $D \in \cK$,
and $C_1 \cap C_2 \subseteq B$ and for every relation symbol $R$ in $\tau$,
if $\bar a\subseteq C_1\cup B\cup C_2$ and $R(\bar a)$ holds then
$\bar a \subseteq C_1 \cup B$ or $\bar a \subseteq C_2 \cup B$(possibly both).

When $D$ is clear from the context we may omit it.
\end{definition}

\section{Abstract Closure Context}
Here we are inside the 0-1-context but without the $\leq_i$ and $\leq_s$
as defined in \S1,  however $\cl^k$ is given. The main result is a sufficient
condition for having 0-1 law or at least convergence.
We have here some amount of freedom, so we give two variants of the main
result of this section: \ref{2.6}, \ref{2.14}, we shall use \ref{2.14}. Thus
on a  reading one may skip Definitions \ref{2.7new} (``possible''),
\ref{2.5} and \ref{2.7A}, Remark \ref{2.5A} and Lemma \ref{2.6} in favour of
the alternative development in Definitions \ref{2.5B}, \ref{2.12new}
and \ref{2.14}. Lemma \ref{2.6A} is needed in both cases and we have made the
two independent at the price of some repetition. We want to ``eliminate
quantifiers'' in a restricted sense: in the simple form we quantify only on
the closure so each $\varphi(\bar x)$ is equivalent to some $\psi_\varphi$ in
which quantifiers are over $\cl^{k,m}(\bar x)$; all this is for a random enough
model where $\cl^{k,m}$ is ``small", still it is not necessarily ``tiny''. The
closure does not need to be in $\cK_\infty$ (though in our application it
is). The quantifier elimination result generalizes the result of
\cite{BlSh:528}. The chief additional ingredient in the proof here is the use
of the addition ($=$ Feferman--Vaught) theorem to analyze a pair of models in
stable amalgamation; this is necessary as we do not have {\em a priori}  bound
on the size of the closure, whereas there we have. Moreover, the argument in
\cite{BlSh:528} is simpler because $<_i$ is defined concretely from a
dimension function and moreover it deal with the ``nice'' rather than almost
nice case.

Note that the ``simply$^*$'' version (\ref{2.20} -- \ref{2.25}) is
used in \S7.

Note that in this section we can forget about the probability distribution:
just deal
with elimination of quantifiers.
Note that the assumption ``$\cl$ is f.o. definable'' (\ref{2.1} clause (d)) is
not serious: if it fails all we have to do is to allow ``$y\in \cl^{k}
(\bar x)$'' as atomic formulas in $\psi_\varphi$.

\begin{context}
\label{2.0}
In this context in addition to $\GK$ (defined in \ref{1.0}) we have an
additional basic operation $\cl$  which is a closure operation for $\cK$ (see
\ref{2.1}), so $\cl$ is in general not defined by Definition \ref{1.3} and
$\leq_i$, $\leq_s$, $\leq_a$ are defined by Definition \ref{2.3} and in
general are not the ones defined in Definition \ref{1.2}. However, we use
$\cK_\infty$ (from \ref{1.2}(1)). Lastly $\nonforkin{}{}_{}^{}$ is as in
\ref{1.16} (can be axiomatized too and moreover generalize to the case of
non--uniqueness, as in \cite{Sh:550}). For simplicity assume $\tau_{\cK}$ (the
vocabulary of $\cK$) is finite with no function symbols. In later
sections (\S4 -- \S7 but not \S3)
saying $\GK$ means $\cl$ is from \S1.
\end{context}

\begin{definition}
\label{2.1}
1) We say $\cl$ is a {\em closure operation} for $\cK$ if for $M\in{\cK}$ and
$k\in\bbn$ the operation $\cl^k (X, M)$ is defined if and only if $X\subseteq
M$ and the operation satisfies:
\begin{enumerate}
\item[(a)] $X\subseteq\cl^k(X, M)\subseteq M$, and $X\subseteq Y\subseteq M
\Rightarrow\cl^k(X,M)\subseteq \cl^k(Y,M)\subseteq M$,
\item[(b)] (i) if $\cl^k(X,M)\subseteq N\subseteq M$ then
$\cl^k(X,N)=\cl^k(X,M)$,\\

(ii) If $X \subseteq N \subseteq M$ then $\cl^k(X,N) \subseteq \cl^k(X,M)$
\item[(c)] for $k\leq\ell$, $\cl^k(X,M)\subseteq \cl^\ell(X,M)$.
\item[(d)] the relation ``$b \in \cl^k(A,M)$'' is preserved by isomorphism.
\end{enumerate}
2) We say that the closure operation $\cl$ is {\em f.o. definable} if (e)
below is true (and we assume this when not said otherwise)
\begin{enumerate}
 \item[(e)] the assertion ``$b\in\cl^k(\{a_0,\ldots a_{\ell-1}\},M)$'' is f.o.
definable in ${\cK}$ that is there is a formula $\psi(y,x_0,...,x_{l-1})$
such that if $M \in \cK$ and $b,a_0,...,a_{l-1} \in M$ then $b \in
cl^k(\{a_0,...,a_{l-1} \},M)$ iff $M \models \psi[y,x_0,...,x_{l-1}]$.
\end{enumerate}
3) We say $\cl$ is {\em transitive} if for every $k$ for some $m$, for every
$X\subseteq M\in\cK$ we have $\cl^k(\cl^k(X,M),M)\subseteq\cl^m(X, M)$.
\end{definition}

\begin{definition}
\label{2.2}
\begin{enumerate}
\item For $X\subseteq M$ and $k,m\in\bbn$ we define $\cl^{k,m}(X,M)$ by
induction on $m$:
\begin{quotation}
$\cl^{k,0}(X,M)=X$

$\cl^{k,1}(X,M)=\cl^{k}(X,M)$

$cl^{k,m+1}(X,M)=cl^{k,1}(cl^{k,m}(X,M),M)$
\end{quotation}
(if we write $\cl^{k,m-1}(X,M)$ and $m=0$ we mean $\cl^{k,0}(X,M)=X$)
\item We say the closure operation $\cl^k$ is $(\ell,r)$-local when:

for $M\in\cK$, $X\subseteq M$ and $Z\subseteq M$ if $Z\subseteq\cl^k(X,M)$,
$|Z|\leq\ell$ then for some $Y$ we have $Z\subseteq Y$, $|Y|\leq r$ and
$\cl^k(Y \cap X, M\restriction Y)=Y$.
\item We say the closure operation $\cl$ is local {\em if} for every $k$, for
some $r$, $\cl^{k}$ is $(1,r)$-local. We say that $\cl$ is simply local if
$\cl^k$ is $(1,k)$-local for every $k$.
\end{enumerate}
\end{definition}

\begin{remark}
\begin{enumerate}
\item Concerning ``possible in $\GK$''(from Definition \ref{2.7new} below),
in the main case $\cM^0_{n,\bar p}$, it
is degenerate, i.e. if $\bar a\subseteq N\in\cK_\infty$, $B\subseteq
N$ then $(N, B,\bar a,k,m)$
is possible. But for the case with the successor relation it has a real role.
\item Note: {\em if} $\cl^k$ is $(1,r)$-local and ``$y\in\cl^k(\{x_1,\ldots,
x_r\},M)$'' is f.o. definable {\em then} for every $m$, $s$ we have ``$y\in
\cl^{k,m}(\{x_1,\ldots,x_s\},M)$'' is f.o. definable.
\item Clearly $\cl^{k,m_1}(\cl^{k,m_2}(X,M))=\cl^{k,m_1+m_2}(X,M)$
and $k_1 \leq k_2 \wedge m_1 \leq m_2 \Rightarrow \cl^{k_1,m_1}(A,M)
\subseteq \cl^{k_2,m_2}(A,M)$.
\item Note that if $\cl^k$ is $(\ell_1,r_1)$-local and $r_2 \geq mr_1$ and
$\ell_2 \leq m\ell_1$ then $\cl^k$ is $(\ell_2,r_2)$--local.
\end{enumerate}
\end{remark}

\begin{definition} [For our 0-1 context $(\cK,cl)$ with $\cl$ as a basic
operation]
\label{2.3}
\
\begin{enumerate}
\item $A\le_i B$ if and only if $A\subseteq B\in {\cK_\infty}$ and for some
$k,m\in\bbn$ and every random enough $\cM_n$ and for every embedding
$g:B\hra\cM_n$ we have $g(B)\subseteq\cl^{k,m}(f(A),\cM_n)$.
\item $A<_s B$ if and only if $A\subseteq B\in {\cK_\infty}$ and for
every $k,m\in
\bbn$ and random enough $\cM_n$ and $f:A\hra M_n$ there is $g$ such that
$f\subseteq g$, and $g:B\hra\cM_n$ with $g(B)\cap\cl^{k,m}(f(A),\cM_n)=f(A)$.
We define $\leq_{pr}$, $\leq^s_m$,$\leq^i_{k,m}$ as in $1.4(5)$, $(7)$, $(8)$
respectively and $A<_a B$ means $A<B$ and $\neg(A<_s B)$.

\item $(\GK,\cl)$ is {\em weakly nice} if for every $A\subseteq C\in
\cK_\infty$, for some $B$ we have $A\leq_i B \leq_s C$.
\item We say $\GK$ (more exactly $(\GK,\cl)$) is {\em
smooth}\footnote{ Smoothness is not used in \cite{Sh:550}, but the
closure there has {\em a priori} bound, so the definitions there will
be problematic here. See more in \cite{Sh:F192}.} when:
\begin{quotation}
\noindent {\em if} $A\subseteq B\subseteq N\in\cK_\infty$, $A\subseteq C
\subseteq N$, $\nonforkin{B}{C}_{A}^{N}$,

\noindent {\em then} $B<_i B\cup C\Leftrightarrow A<_i C$
\end{quotation}
(note that $\Leftarrow$ is always true).
\item We say that $\cl^k$ is $r$-transparent if
$$
A\leq_i B\ \&\ |B|\leq r\qquad\Rightarrow\qquad \cl^k(A,B)=B.
$$
We say that $\cl$ is transparent if for every $r$ for some $k$ we
have: $\cl^k$ is $r$-transparent.
We say that $\cl$ is simply transparent if for every $k$,  $\cl^k$ is
$k$-transparent.
\end{enumerate}
\end{definition}

\begin{fact}
\label{2.4}
Assume $\GK$ is a $0-1$ context (see \ref{1.0}) and $\cl$ is defined in
\ref{1.3} {\em  then }

\begin{enumerate}
\item[$(\alpha)$] $\cl$ is a closure operation for $\cK_{\infty}$ (see
Def.\ref{1.1}(1)),
\item[$(\beta)$] $\cl$ is f.o. definable (for $\cK$),
\item[$(\gamma)$] $\cl^{k,m}$ as defined in \ref{1.3}(c) and as defined
\ref{2.2}
are equal,
\item[$(\delta)$] $\cl$ is transitive,
\item[$(\vep)$] $\cl$ is simply local (see Def.\ref{2.2}(2),(3)),
\item[$(\zeta)$] $\cl$ is transparent, in fact $\cl^k$ is $k$-transparent for
every $k$,
\item[$(\eta)$] $\leq_i$ as defined in \ref{2.3}(1) and in \ref{1.2} are equal,
\item[$(\theta)$] If in $\S 1$, $\GK$ is weakly nice (see Def.\ref{1.5}) then
$(\cK_{\infty},\cl)$ is weakly nice by Def.\ref{2.3}(3); if so then
$\leq_s$ as
defined in \ref{2.3}(2) and \ref{1.2}(4) are the same and $<_a$ as defined in
\ref{2.3}(2) and in \ref{1.2}(6) are equal.
\end{enumerate}
\end{fact}

\noindent PR
OOF.\quad
$(\alpha)$\qquad We have to show that $(\cK_\infty,\cl)$ from $\S 1$
satisfies clauses $(a),(b),(c),(d)$ from Def. \ref{2.1}(1).\\
$(a)$ By the Def. \ref{1.3} of $(\cK_\infty,\cl)$ the
following
holds: trivially
$X \subseteq \cl^k(X,M) \subseteq M$. Assume $X \subseteq Y \subseteq M$.
If $b \in \cl^k(X,M)$ then
for some $B$,
$|B| \leq k$ and $b \in B$, $X \cap B \leq_i B$ by Def.$1.6$. As $X \subseteq
Y$
and $X \cap B \leq_i B$ we obtain $Y \cap B \leq_i B$
by
Fact \ref{1.6}(4). So $B\subseteq \cl^k(Y,M)$ witnessing that
$b\in \cl^k(Y,M)$. Hence $cl^k(X,M) \subseteq cl^k(Y,M)$.\\

$(b)$ $(i)$ First, let's show $\cl^k(X,N)\subseteq \cl^k(X,M)$.
If $b\in \cl^k(X,N)$ then let $B$ witness it and we have $b\in B$,
$B\subseteq N$, $B\cap
X\leq_i B$, $|B|\leq k$. As $N\subseteq M$ the witness $B$ is in $M$, $B\cap
X \leq_i B$ so $b\in \cl^k(X,M)$. Second we will show that $\cl^k(X,M)
\subseteq
\cl^k(X,N)$: if $b\in \cl^k(X,M)$ then there is $B$ witnessing it such that
$b\in B\subseteq M$,$B\cap X \leq_i B$, $|B|\leq k$. Now clearly $B\subseteq
\cl^k(X,M)$ hence by assumption $B\subseteq N$  so
$b\in B \subseteq N$, $B \cap X\leq_i B$, $|B|\leq k$ and so $B$ witnesses
$b \in \cl^k(X,N)$. So we get the result.\\

$(ii)$ Included in the proof of clause $(i)$.

$(c)$ It follows immediately that $(\cK,\cl)$ holds by Def.$1.6$.

$(d)$ Easy.

\noindent $(\beta)$\qquad We show that $(\cK,\cl)$ is f.o. definable. By
Def.\ref{2.1}(d) this
means that for each $\ell$, there is a formula $\psi(y,x_0,...,x_{\ell-1})$
such that if
$M\in \cK$ and $b,a_0,...a_{\ell-1} \in M$ then: $b \in
\cl^k(\{a_0,...,a_{\ell-1} \}, M)$
iff $M\models \psi(b,...,a_{\ell-1})$. It suffice to restrict ourselves to
the case $\langle b_0,...,b_{\ell-1} \rangle$ is with no repetition.

Let $ {\mathfrak B}=\{(B,\bar b): B \in \cK_\infty$ has $\leq \ell+1$
elements, $\bar b$
is a sequence of length $\leq \ell+1$ listing the elements of $B$
without repetitions$\}$. On ${\mathfrak B}$ the relation $\cong$ (isom.) is
defined. We say $(B',\bar b') \cong (B{''},\bar b{''})$ if there is an
isomorphism $h$ from $B'$ onto $B''$ mapping $\bar b'$ onto $\bar b''$. Now
$\cong$ is an equivalence relation on ${\mathfrak B}$. ${\mathfrak B}/ \cong$
is finite. So let $\{(B_i,\bar b_i): i<i^* \}$ be a set of representatives.
Now $i^*$ is finite as $\tau$ is finite
(actually locally finite suffice). Let when $k=k_i=|B_i|=lg(\bar b_i)$
\begin{multline*} \varphi_i(x_0,...,x_{k})=\\
\bigwedge \{ \theta(x_0,...,x_{k}):
\theta \mbox{ is a basic formula (possibly with dummy variables) and}\\
 B_i \models
\theta [b_0,...,b_{k}] \}.
\end{multline*}
Lastly
\[\begin{array}{l}
\psi(y,x_0,...,x_{\ell-1})=\bigvee\limits_{m<\ell}y=x_m\vee \\
\qquad\bigvee\{(\exists z_0,\ldots,z_{k-1})
(\bigwedge\limits_{m<\ell}\bigvee\limits_{t<k}x_m=z_t \wedge y=z_k \wedge
\varphi_i(y,z_0,...,z_{k-1})):\\
\qquad\qquad B_i\mbox{ has exactly $k+1$ members and }B_i\upharpoonright
\{b^i_t: t<\ell\} \leq_i B_i \} 
  \end{array}\]

\noindent $(\gamma)$\qquad Trivial.

\noindent $(\delta)$\qquad By \ref{1.11}.

\noindent $(\vep)$\qquad Now, we will show that $(\cK_\infty,\cl)$ is simply
local.  For this we have to show that $\cl^k$ is $(1,k)$-local for every
$k$. Let $X \subseteq M \in \cK_\infty$ be given and $Z\subseteq \cl^k(X,M)$
such that $|Z| \leq 1$. If $Z=\emptyset$ let $Y=\emptyset$. So assume
$Z=\{y\}$. As $y \in Z \subseteq \cl^k(X,M)$ there is a witness set $Y$ for
$y \in \cl^k(X,M)$ so $Y \cap X \leq_i Y$, $|Y| \leq k$. As $Y \cap X \leq_i
Y$, clearly $\cl(X \cap Y,Y)=Y$ and $Z=\{y \}\subseteq Y$ and $|Y| \leq k$
so we are done. 

\noindent $(\zeta)$\qquad Trivial by the definition of $\cl$ (Def.$1.6$) and
of transparency (Def.\ref{2.3}(5)).

\noindent $(\eta)$\qquad First assume $A \leq_i B$ by Def. \ref{2.3} and we
shall prove that $A \leq_i B$ by Def. \ref{1.2}. So for some $k,m$ we have:
\begin{enumerate}
\item[$(*)$] for every random enough $\cM_n$ and embedding $g: B \hra \cM_n$
we have $g(B) \subseteq \cl^{k,m}(g(A),\cM_n)$.
\end{enumerate}
Let $\vep > 0$. Let $\cM_n$ be random enough and $f:A \hra \cM_n$.
By (*) and \ref{1.11} if $g$ is an embedding 
of $B$ into $\cM_n$ extending $f$
then we have $g(B) \subseteq \cl^{k^m}(g(A),\cM_n)$, hence\\
$|\ex(f,B,\cM_n)| \leq |\cl^{k^m}(g(A),\cM_n)|^{|B\setminus A|}$.
Let $\zeta=\vep/(|B\setminus A|+1)$, now if $\cM_n$ is random enough,
then by \ref{1.12} for every
$g: B \hra \cM_n$ we have $|\cl^{k^m}(g(A),\cM_n)|^{|B\setminus A|}
\leq n^{\zeta}$, hence $\ex(f,B,\cM_n)| \leq |n^{\zeta}|^{|B\setminus A|}
\leq n^{\vep}$. As $\vep>0$ was arbitrary, we have proved that $A \leq_i B$
by Def.\ref{1.2}.

Next assume $A \leq_i B$ by Def. \ref{1.2} and we shall prove that $A \leq_i B$
by Def.\ref{2.3}. Choose $k=|B|$ and $m=1$, so $\cl^{k,m}=\cl^k$. So let
$\cM_n$ be random enough, and $g: B \hra \cM_n$. Recall that
$\cl^k(g(A),\cM_n)=\cup \{C:C\subseteq \cM_n, |C| \leq k$ and $C \cap A
\leq_i C \}$, so $g(B)$ can serve such $C$, hence $g(B) \subseteq
\cl^k(g(A),\cM_n)$.

\noindent $(\theta)$\qquad We shall use clause $(\eta)$ freely.
First assume that $\cK$ is weakly nice by Def.\ref{1.5} and we shall prove
that $(\cK,\cl)$ is weakly nice by Def.\ref{2.3}(3). So assume $A \leq
B$. We can find $C$ such that $A \leq_i C \leq B$ and for no $C'$, $A \leq_i
C' \leq B$, $C \subset C'$; exist as $A \leq_i A \leq B$ and $B$ is
finite. By \ref{1.6}(2) for no $C'$, do we have $C <_i C' \leq B$ hence $C
\leq_s B$ by Def.\ref{1.2}, so it is enough to prove that $C \leq_s B$ by
Def.\ref{2.3}(2), and w.l.o.g. $C \ne B$ so $C <_s B$. Let $k,m$ be
given. As we are assuming  that $\cK$ is weakly nice by Def.\ref{1.5} and $C
<_s B$ by Def.\ref{1.2}(4) we have that there is an $\vep \in \bbr^+$ such
that 
$$
1=\lim\limits_n\left[\prob\left(
\begin{array}{l}
\mbox{if } f_0:A \hra \cM_n\ \mbox{ then there is } F \mbox{ with }
|F|\ge n^\varepsilon\ \mbox{ and }\\
\mbox{ (i) } f\in F\imply f_0 \subseteq f: B \hra \cM_n\\
\mbox{ (ii) } f'\neq f''\in F\imply\Rang(f')\cap\Rang(f'')=\Rang (f_0)
\end{array}
\right)\right].
$$
As $\cM_n$ is random enough and $f:A \hra \cM_n$, there is $F$ as
above for $B$ with $|F| \geq n^{\vep}$; but by \ref{1.11} also
\[|\cl^{k,m}(f(A),\cM_n)| \leq |\cl^{\ell}(f(A),\cM_n)|\]
for $l=k^m$ and by \ref{1.4} we have
\[|\cl^{{k^m}}(f(A),\cM_n)|< n^{\vep}\]
so $|\cl^{k,m}(f(A),\cM_n)|< n^{\vep}$.

As the sequence $\langle\Rang (g\setminus\Rang (f):g\in F\rangle$
list a family of $\geq n^{\vep}>|\cl^{k,m}(f(A),\cM_n)|$  pairwise disjoint
subsets of $\cM_n$, for some $g \in F$, we have: $\Rang(g)\cap\Rang(f)$
is disjoint to $\cl^{k,m}(f(A),\cM_n)$. So $g$ is as required in
Def.\ref{2.3}(2); so we have finished by  proving $C \leq_s B$ by
Def.\ref{2.3}, hence we have finished proving $(\cK_\infty,\cl)$ is weakly
nice according to Def.\ref{2.3}(3).\\ 
So we have proved the implication between the two version of weakly nice.
Second, assuming $\GK$ is weakly nice by Def.\ref{1.5}, we still have to
say why the two version of $\leq_s$ (by Def.\ref{1.2}(4) and by \ref{2.3}(2))
are equivalent. Now if $C \leq_s B$ by Def.\ref{1.2}(4) then
$C \leq_s B$ by Def.\ref{2.3}(2) has been proved inside the proof above that
$\cK$ weakly nice; by Def.\ref{1.2}(3) implies $(\cK,\cl)$ is weakly
nice by Def.\ref{2.3}(3). Lastly assume $A \leq_s B$ by Def.\ref{2.3}(2),
now if $A <_i C \leq B$ we get a contradiction
directly from Def.\ref{2.3}(2): but this confirm $A \leq_s B$ according to
Def.\ref{1.2}(4).

Lastly we leave the statement on $<_a$ to the reader.
\hfill\rqed$_{\ref{2.4}}$

\begin{remark}
\label{2.7tikon}
\begin{enumerate}
\item Note that the assumption ``$\GK$ is weakly nice'' is very natural
in the applications we have in mind.
\item Why have we not prove the equivalence of the two versions of
weakly nice in \ref{2.4}$(\theta)$? We can define the following $0-1$
context $\GK$: let $\cM_n$ be $\cM^0_{n,\bar p}$ if $n$ is even with
$p_n=1/n^\alpha$, $\alpha \in (0,1)_{\bbr}$ irrational (except
$p_1=1/2^{\alpha})$ and $\cM_n$ is the random graph with probability
$1/2$ if $n$ is odd. Now in \S1, $\cK_{\infty}$ is the family of finite
graphs, and $A \leq_i B$ iff $A=B$ (using the odd $n$-s). Hence
$\cl^k(A,M)=A$ so clearly $A<B \Rightarrow A<_s B$ according to \ref{1.2}
hence weakly niceness by \ref{2.3}(3) holds trivially but weakly niceness by
Def.\ref{1.5} fail.
\item Note that in Definitions \ref{2.7new}, \ref{2.5}, \ref{2.5B} below the
``universal''
demand speak on a given situation in random enough $\cM_n$ whereas
the ``existential demand'' implicit in goodness deal with extensions
of an embedding into $\cM_n$.
\item We would like to show that for every formula $\varphi (\bar x)$
(f.o. in the vocabulary $\tau_{\cK}$) there are (f.o.) $\psi_{\varphi}(
\bar x)$ and  $k=k_{\varphi}$, $m=m_{\varphi}$ such that\\
\noindent $(*)_{\varphi}$ for every random enough $\cM_n$ and $\bar a
\in\ ^{lg(\bar x)}\cM_n$ we have
$\cK \models \varphi[\bar a] \Leftrightarrow \cM_n \restriction
cl^{k,m}(\bar a,\cM_n) \models \psi_{\varphi}(\bar x)$.\\
\noindent Naturally enough we shall do it by induction on the
quantifier depth of $\varphi$ and the non-trivial case is
$\varphi(\bar x)=(\exists y)
\varphi_1(\bar x,y)$, and we assume $\psi_{\varphi_1}(\bar x,y)$,
$k_{\varphi_1}$, $m_{\varphi_1}$ are well defined. So we
should analyze the situation: $\cM_n$ is random enough, $\bar a \in\
{}^{lg(\bar x)}(\cM_n)$, $\cM_n \models \varphi[\bar a]$ so there is $b \in
\cM_n $ such that $\cM_n \models \varphi_1[\bar a,b]$, and we split it to
two cases according to the satisfaction of a suitable statement on a
suitable neighbourhood of $\bar a$ i.e., $cl^{k',m'}(\bar a, \cM_n)$. If $b$
belongs to a small enough neighbourhood of $\bar a$ this should be clear. If
not we would like to find a suitable situation (really a set of possible
situation, with a bound on their number depending just on $\varphi$) to
guarantee the existence of an element $b$ with
$\cl^{k_{\varphi_1},m_{\varphi_1}}(\bar a b,\cM_n)$ satisfying
$\psi_{\varphi_1}(\bar a,b)$. Now in general the 
$\cl^{k_{\varphi_1},m_{\varphi_1}}$ can
be of large cardinality (for $\varphi$, i.e. depending on $\cM_n$). In the
nice case we are analyzing, to find such a witness $b$ outside a small
neighbourhood of $\bar a$ it will suffice to look at
$\cl^{k_\varphi,m_\varphi}(\bar ab,\cM_n)$ essentially with small
cardinality. Why only essentially? As may be
$\cl^{k_{\varphi_1},m_{\varphi_2}}(\bar a,\cM_n)$ is already large,
so what we should have is something is like:
$\cl^{k_{\varphi_1},m_{\varphi_1}}(\bar ab,\M_n) \setminus
\cl^{k_{\varphi_1,m_2}}(\bar a,\cM_n)$
can be replace by a set of small cardinality. For this we need
$\nonforkin{}{}_{}^{}$
(the relation free amalgamation) to hold, possibly replacing
$\cl^{k_{\varphi_1},m_2}(\bar a,\cM_n)$ by a subset (in $\S3$ we can make it
arbitrary, here quite definable)
and the amalgamation base has an a priori bound.
By the addition theorem we may replace $(B^*,b)_{b \in B}$
by similar enough
$(B',b)_{b \in B}$ (in particular when $B^* \in \cK_{\infty}$ so we need to
express in such situation something like $B^*$ exists over $B$
(we can say such $B$ exists by clause $(b)$ of \ref{2.7new}(4) using
quantifiers on $\cl^{k,m}(\bar a, \cM_n))$. Well, $B \leq_s B^*$
is good
approximation. But this does not say that $\cl(\bar ab,\cM_n)$ is suitable.
So we need to say first that the closure of $\bar a b$ in essentially
$B^* \cup B_2$ where $B_2=cl^{k_{\varphi_1},m_2}(\bar a,\cM_n)$, obeys
a version of the addition theorem, and secondly that $B^*$ sit in $\cM_n$
in a way where the closure is right.
All this is carried in Def.\ref{2.7new}(4) (of good saying: we have a tuple
in a
situation which exist whenever a copy of $B$ as above exist) and \ref{2.5}
(when there are $B$ etc. as above). The proof is carried in \ref{2.6}.
\item Defining good, by demanding the existence of the embedding
$g:B^* \hra\cM_n$
extending $f:B\hra\cM_n$, we demand on $f$ only little: it is an
embedding. We may impose requirements of the form
$cl^{k_i,m_i}(f(B_i),\cM_n) \subseteq f(B)$ or $cl^{k_i,m_i}(f(B_i),\cM_n)
\cap f(B)=f(C_i)$ for some $B_i$,$C_i \subseteq B$. This make it easier for
a tuple to be good. Thus giving a version of almost nice covering more cases.
In other possible strengthening we do not replace $B^{*}$ by
$B' \in \cK_{\infty}$ of bounded cardinality but look at it as a family of
possible ones all similar in the relevant sense. On the other hand we may
like simpler version which are pursued in \ref{2.12new}, \ref{2.14}.
\item Note that if $\cl^k$ is $r$-transparent and $A \subseteq M
\in \cK$ then $\cl^k(A,M) \supseteq \cup \{ C \subseteq M: C \cap A \leq_i C$
and $|C| \leq r \}$. [Why? if $C \subseteq M$, $C \cap A \leq_i C$ and
$|C| \leq r$ then : first $\cl^k(C \cap A,C)=C$ as $\cl^k$ is
$r$-transparent;
second $\cl^k(C \cap A,C) \subseteq \cl^k(C \cap A,M)$ by $(b)(ii)$ of Def.
 \ref{2.1}
(1), third $\cl^k(C \cap A,M) \subseteq \cl^k(A,M)$ as $C \cap A
\subseteq M$ by clause $(a)$ of Def.\ref{2.1}(1); together we are done]. Note
that if $\cl^k$ is $(1,r)$-local we can prove the other inclusion. So
obviously if $(\GK,\cl)$ is simply local and simply transparent
(and  $\tau_{\mathcal K } $ is finite or at least locally finite of course),
then $\cl$ is f.o. definable. If we omit the simple we can eliminate the 
assumption $\cl$ is f.o. definable in \ref{2.6}, \ref{2.14}.
\end{enumerate}
\end{remark}

\begin{definition}
\label{2.7new}
\begin{enumerate}
\item We say $(N,B,\bar B,k)$ is possible for $(\GK,\cl)$ if:
\begin{enumerate}
\item[(a)] $\bar B=\langle B_i: i< lg(\bar B)\rangle$, $B_i\subseteq
N\in\cK_\infty$, $B\subseteq N$ and $\cl^k(B_i,N)
\subseteq B_{i+1}$ for $i<\lg(\bar B)-1$
\item[(b)] it is not true that:

for every random enough $\cM_n$, for no embedding $f: N\hra\cM_n$, do we have:
$$
\mbox{for }i<\lg(\bar B)-1,\ \cl^k(f(B_i),\cM_n)\subseteq f(\cl^k(B_i),
N)\cup\cl^{k}(f(B),\cM_n).
$$
\end{enumerate}
\item If we write $(N,C,B,k)$ we mean $(N,C,\langle B,\cl^{k}(B,N)\rangle,k)$.
\item We say $(N,B,\bar a,k,m)$ is possible for $\GK$ if $(N,B,\bar B,k)$ is
possible for $\GK$ where $\bar B=\langle\cl^{k,i}(\bar a,N):i\leq m\rangle$.
\item We say that the tuple $(B^*,B,B_0,B_1,k,m_1,m_2)$ is good for
$(\cK,\cl)$
if
\begin{enumerate}
\item[(a)] $B\leq B^* \in \cK_{\infty}$ and, $B_0\leq B_1 \leq B^* \in
\cK_{\infty}$
\item[(b)]
for every random enough $\cM_n$ we have:
if $f: B\hra\cM_n$ then there is an extension
$g$ of $f$ satisfying $g:B^*:\hra\cM_n$ and
\end{enumerate}

\begin{enumerate}
\item[($\alpha$)] $g (B^*)\cap\cl^{k,m_2}(f(B),\cM_n)=f(B)$,
\item[($\beta$)] $\cl^{k,m_1}(g(B_0),\cM_n)\subseteq g(B_1)\cup
\cl^{k,m_2}(g(B),\cM_n)$
\item[($\gamma$)] $\nonforkin{\cM_n\restriction g(B^*)}{\cM_n
\restriction\cl^{k,m_2}(f(B),\cM_n)}_{\cM_n\restriction f(B)}^{\cM_n}$,
\end{enumerate}
\end{enumerate}
\end{definition}

\begin{definition}
\label{2.5}
The 0-1 context $\GK$ with closure $\cl$ (or the pair $(\GK,\cl)$ or $\GK$
when $\cl$ is understood) is {\em almost nice} if it is weakly nice and 
\begin{enumerate}
\item[(A)] {\em the universal demand:}

for every $k,m_0$ and $\ell,\ell'$ there are
$$
\hspace{-0.5cm}m^\ast = m^\ast(k,m_0,\ell,\ell')>m_0,\ k^*=k^*(k,m_0,\ell,
\ell')\geq k\ \mbox{ and }t=t(k,m_0,\ell,\ell')
$$
such that, for every random enough $\cM_n$ we have:

{\em if} $\bar a\in {}^\ell |\cM_n|$ and $b\in\cM_n\setminus\cl^{k^\ast,
m^\ast}(\bar a,\cM_n)$

{\em then} there are $m_2\in [m_0,m^*]$ and $m_1\le m^\ast-m_2$ and
$B\subseteq\cl^{k,m_1}(\bar a,\cM_n)$ and $B^*\subseteq\cM_n$ such that:
\item[($\alpha$)] $|B|\le t$ and $\bar a\subseteq B$,
\item[($\beta$)]  $B^*=[\cl^{k,m_0}(\bar ab,\cM_n)\setminus\cl^{k,m_2}(B,
\cM_n)]\cup B$ so necessarily $b\in B^*$ and $\bar a\subseteq B^*$, (see
\ref{2.5A} below)
\item[($\gamma$)] $B<_s B^*$ or at least:\\
for every first order formula $\varphi=\varphi(\ldots,x_a,\ldots)_{a\in B}$ of
quantifier depth $\leq\ell'$ there is $B'$ such that $B<_s B'$ (so
$B'\in\cK_\infty$) and
$$
B^*\models\varphi(\ldots,a,\ldots)_{a\in B}\quad\mbox{ iff }\quad B'\models
\varphi(\ldots,a,\ldots)_{a\in B},
$$
\item[($\delta$)] $\nonforkin {\cM_n\upharpoonright B^*}{\cM_n\upharpoonright
\cl^{k,m_2}(B,\cM_n)}_{\cM_n\upharpoonright B}^{\cM_n}$,
\item[$(\varepsilon)$] $(B^*,B,\bar ab,B^*\cap\cl^{k,m_0}(\bar ab,
\cM_n)),k,m_0,m_2)$ is good for $(\cK,cl)$ or at least for some $B'$, $B''$
we have\footnote{$M_1 \equiv_{\ell'} M_2$ means;
$M_1,M_2$ satisfy
the same f.o. sentences of quantifier depth $\leq \ell'$}
\end{enumerate}
\begin{enumerate}
\item[(i)] $(B',B,\bar ab,B'',k,m_1,m_2)$ is good for $(\cK,cl)$
\item[(ii)] $(B^*,\cl^{k,m_0}(\bar ab,\cM_n)\cap B^*), b,c)_{c\in B}
\equiv_{\ell'}(B',B'',b,c)_{c\in B}$,
\item[$(\zeta)$] for $m\leq m_{\varphi_1}$ we have
$$
\cl^{k_{\varphi_1},m}(\bar ab,B^*)=B^*\cap\cl^{k_{\varphi_1},m}(\bar ab,
\cM_n).
$$
\end{enumerate}
\end{definition}

\begin{definition}
\label{2.7A}
If in Def \ref{2.5} above, $k^*=k$ in clause ($A$)
{\em then } we add ``$k$--preserving''.
\end{definition}

\begin{remark}
\label{2.5A}
\begin{enumerate}
\item Note that if $\cK=\cK_\infty$ and $\cl$ is local (or just $\cl^k$ is
$(l_k, r_k)$-local for each $k$) (which holds in the cases we are interested
in) then in clauses $(\gamma)$, $(\varepsilon)$ of (A) in Def.\ref{2.5}
above the two possibilities are close. 
\item Why in \ref{2.5}(A)$(\beta)$ we have ``necessarily $b\in B^*$''?
Because
$$
b\in\Rang(\bar ab)\subseteq\cl^{k,m_0}(\bar ab,\cM_n)\quad\mbox{ and}
$$
$$
\begin{array}{ll}
\cl^{k,m_2}(B,\cM_n)&\subseteq\cl^{k,m_2}(\cl^{k,m_1}(\bar a,\cM_n))\subseteq
\cl^{k,m_1+m_2}(\bar a,\cM_n)\\
\ &\subseteq\cl^{k,m^*}(\bar a,\cM_n)\subseteq\cl^{k^*,m^*}(\bar a,\cM_n)\\
\end{array}
$$
and $b$ does not belong to the later.
\item Why do we use $\cl^{k,m_2}(B,\cM_n)$? Part of our needs is that this set
is definable from $B$ without $b$.
\item In clause $(\gamma)$, Definition \ref{2.5} clause (A), there is one $B'$
for all such $\varphi$ (Why? As the set of f.o. formulas of quantifier
depth $\ell$
is closed under Boolean combinations) so for some $B'\in\cK_\infty$ we have
$B\leq_s B'$, and $(B',c)_{c\in B}\equiv_\ell(B^*,c)_{c\in B}$. So we could
have phased clause $(ii)$ of $(A)(\vep)$ in the same way as clause
$(\gamma)$.
\end{enumerate}
\end{remark}

In our main case, also the following variant of the property applies (see
\ref{2.20new} below).

\begin{definition}
\label{2.5B}
1) We say that the quadruple $(N,B,\langle B_0,B_1\rangle,k)$ is {\em
simply good} for $(\GK,\cl)$ if
($B$, $B_0$, $B_1\leq N\in \cK_\infty$ and)
for every random enough $\cM_n$, for every embedding $f: B\hra\cM_n$ there is
an extension $g$ of $f$ satisfying $g: N\hra\cM_n$ such that:
\begin{enumerate}
\item[(i)]   $g(N)\cap\cl^k(f(B),\cM_n)=f(B)$,
\item[(ii)]  $\nonforkin{g(N)}{\cl^k(f(B),\cM_n)}_{f(B)}^{}$,
\item[(iii)] $\cl^k(g(B_0), \cM_n)\subseteq g(B_1)\cup \cl^k(g(B), \cM_n)$
\end{enumerate}
(natural but not used is $\cl^k(g(B_0),\cM_n)\cap g(N)=g(\cl^k(B_0,N))$). If
we write $B_0$ instead $\langle B_0, B_1\rangle$, we mean $B_1=N$.

\noindent 2) We say that $(N, B, \langle B_0, B_1\rangle, k, k')$ is {\em
simply good} if part (1) holds replacing (iii) by
\begin{enumerate}
\item[(iii)$'$] $\cl^k(g(B_0), \cM_n) \subseteq g(B_1)\cup
\cl^{k'}(g(B), \cM_n)$.
\end{enumerate}
\end{definition}

\begin{definition}
\label{2.12new}
 1) The 0-1 context with closure $(\GK,\cl)$ is {\em simply almost
nice\/} if it is weakly nice and 
\begin{enumerate}
\item[(A)] {\em the universal demand:}

for every $k$ and $\ell,\ell'$ there are
$$
m^\ast=m^\ast(k,\ell,\ell'),\quad k^*=k^*(k,\ell, \ell')\geq k\ \mbox{ and }t
= t(k,\ell,\ell')
$$
such that for every random enough $\cM_n$ we have:

{\em if} $\bar a\in {}^\ell |\cM_n|$ and $b\in\cM_n\setminus\cl^{k^*,m^*}(\bar
a,\cM_n)$

{\em then} there are $B\subseteq\cl^{k^*,m^*}(\bar a,\cM_n)$ and
$B^*\subseteq\cM_n$ such that:
\begin{enumerate}
\item[$(\alpha)$] $|B|\le t$ and $\bar a\subseteq B$ and $\cl^k(B,\cM_n)
\subseteq\cl^{k^*,m^*}(\bar a,\cM_n)$,
\item[($\beta$)]  $B^*=[\cl^{k}(\bar ab,
\cM_n)\setminus\cl^{k}(B,\cM_n)]\cup B$

(or at least $B^*\supseteq [\cl^k(\bar a b, \cM_n)\setminus \cl^k(B, \cM_n)]
\cup B$),
\item[($\gamma$)] $B<_s B^*$ (so $B^*\in\cK_\infty$) or at least for every
first order formula $\varphi=\varphi(x_b,\ldots,x_a,\ldots)_{a\in B}$ of
quantifier depth $\leq\ell'$ there is $B'$ such that $B<_s B'$ (so $B'\in
\cK_\infty$) and:
$$
B^*\models\varphi(b,\ldots,c,\ldots)_{c\in B}\quad\mbox{ iff }\quad B'\models
\varphi(b,\ldots,c,\ldots)_{c\in B}
$$
(or even , but actually equivalently, $(B^*, b,\ldots,c,\ldots)_{c\in B}
\equiv_{\ell'}(B',b,\ldots,c,\ldots)_{c\in B}$),
\item[($\delta$)]  $\nonforkin{\cM_n\upharpoonright B^*}{\cM_n
\upharpoonright\cl^{k}(B,\cM_n)}_{\cM_n\upharpoonright B}^{\cM_n}$
\item[$(\varepsilon)$] $B^*\in\cK_\infty$ and $(\cM_n\restriction B^*,B,\bar
ab,k)$ is simply good for $(\GK,\cl)$ or at least for some $B'$, $b'$ we have:
\begin{enumerate}
\item[(i)]  $(B',B,\bar ab',k)$ is simply good for $(\GK, \cl)$ and
\item[(ii)] $(B^*,b,\ldots,c,\ldots)_{c\in B}\equiv_{\ell'}(B',b',\ldots,c,
\ldots)_{c\in B}$.
\end{enumerate}
\end{enumerate}
\end{enumerate}
\noindent 2) If above always $k^*=k$ we say: $\GK$ is simply almost nice depth
preserving.

\noindent 3) We say that $(\GK,\cl)$ is {\em simply nice} (i.e. omitting the
almost) if \ref{2.12new}(1) holds but we omit clause $(\varepsilon)$ and add
\begin{enumerate}
\item[(B)] if $B<_s B^*$ and $k\in\bbn$ then
$(B^*,B,B^*,k)$ is simply good.
\item[(C)] $\cK_\infty=\cK$ (or at least if $A \in \cK_\infty$ and
$k,m \in \bbn$
then for any random enough $\cM_n$ for any $f: A \hra \cM_n$, $\cl ^{k,m}(A,
\cM_n) \in \cK_\infty$.
\end{enumerate}
Similarly in Definition \ref{2.5} for ``nice''.
\end{definition}

\begin{remark}
\label{2.12Anew}
1) In \ref{2.12new}(1) we can weaken the demands (and call $(\cK,\cl)$
simply$^{\otimes}$ almost nice): get also
$k^\otimes=k^\otimes(k, \ell, \ell')\in\bbn$ replace in clause $(\beta)$
$\cl^k(B,\cM_n)$ by $\cl^{k^\otimes}(b,\cM_n)$ and
replace $(\varepsilon)$ by
\begin{enumerate}
\item[($\varepsilon'$)] $(B',B,\bar{a}b,k,k^\otimes)$ is simply good for
$(\GK, \cl)$ (see \ref{2.5B}(2)) {\em or } at least for some $B'$, $b'$
we have:
\begin{enumerate}
\item[(i)] $(B', B, \bar ab', k, k^\otimes)$ is simply$^{\otimes}$ good
\item[(ii)] $(B^*, b, \ldots, c, \ldots)_{c\in B} \equiv_{\ell'} (B',
b', \ldots, c, \ldots)_{c\in B}$
\end{enumerate}
\end{enumerate}
The parallel change in \ref{2.12new}(2) (that is defining simply$^{\otimes}$
good) is
\begin{enumerate}
\item[(B)$'$] for every $k, \ell\in \bbn$ for some $k^\otimes= k^\otimes
(k, \ell)\in \bbn$ we have: if $B <_s B^*$ and $|B| \leq \ell$, then
$(B^*, B, B^*, k, k^\otimes)$ is simply good.
\end{enumerate}
This does not change the conclusions i.e (\ref{2.12new}, \ref{2.14},
\ref{2.20new}, \ref{2.7}).

\par \noindent
2) We can change Definition \ref{2.5} as we have changed Definition
\ref{2.12new}(1) in \ref{2.12new}(3) and/or in \ref{2.12Anew}(1).

\par \noindent
3) We can without loss of generality demand in \ref{2.12new}(1)(A) that
$m^*(k, \ell, \ell')=1$ at the expense of increasing $k^*$, as if
$\cl^{k^{**}}(\bar a, M)\supseteq \cl^{k^*, m^*}(\bar a, M)$ whenever $\bar
a\in {}^\ell |M|$, $M\in \cK$ then $k^{**}$ will do.

\par \noindent
4) We can omit clause $(\gamma)$ in Def.\ref{2.12new}(1), but it is natural.
Similarly in Def.\ref{2.5} (i.e. those omitting do not change the
later claims).

\par \noindent 5) In Def.  \ref{2.12new} we can omit $m^*$ if $\cl$ is
transparent by increasing $k^*$ (that is $m^*=1$) 
\end{remark}

Lemma \ref{2.6A} below (the addition theorem, see \cite{CK} or \cite{Gu}
and see more \cite{Sh:463}) is
an immediate a corollary of the well known addition theorem; this is the
point where $\nonforkin{}{}_{}^{}$ is used.

\begin{lemma}
\label{2.6A}
For finite vocabulary $\tau$ and f.o. formula
(in $\tau$) $\psi(\bar z,\bar z^1,\bar z^2)$,
$\bar z=\langle z_1,\ldots,z_s\rangle$, there are $i^*\in\bbn$
and $\tau$- formulas
$\theta^1_i(\bar z, \bar z^1)=\theta^1_{i,\psi}(\bar
z,\bar z^1)$,
$\theta^2_i(\bar z,\bar z^2)=\theta^2_{i,\psi}(\bar z,\bar z^2)$
for $i<i^*$, each of quantifier depth at most that of
$\psi$ such that:
\begin{quotation}
\noindent {\em if} $N$ is $\tau$-model, $\nonforkin{N_1}{N_2}_{N_0}^{N}$,
$N_1\cap N_2= N_0$,
$N_1\cup N_2=N$ and the set of elements of $N_0$ is $\{c_1,\ldots,c_s\}$,
$\bar c=\langle c_1,\ldots,c_s\rangle$ and $\bar c^1\in{}^{\lg\bar z^1}(N_1)$
and $\bar c^2 \in {}^{\lg \bar z^2}(N_2)$

\noindent {\em then:}
$$
\hspace{-1.2cm}N\models\psi[\bar c,\bar c^1,\bar c^2]\quad\mbox{ iff
}\quad\mbox{ for some }i<i^*, N_1\models\theta^1_i[\bar c,\bar c^1]\mbox{ and
} N_2\models\theta^2_i [\bar c,\bar c^2].
$$
\end{quotation}
\end{lemma}

\begin{mainlemma} [Context as above]
\label{2.6}
Assume $(\GK,\cl)$ is almost nice and $\cl$ is f.o. definable.

\noindent 1) Let $\varphi (\bar x)$ be a f.o. formula
in the vocabulary $ \tau_{\mathcal K}$. {\em Then} for some
$m_\varphi\in\bbn$ and $k=k_\varphi\ge\lg(\bar x)+q.d.(\varphi(\bar x))$ and
for some f.o. $\psi_\varphi (\bar x)$ we have:
\begin{enumerate}
\item[$(\ast)_{\varphi}$] for every random enough $\cM_n$ and
$\bar a\in {}^{\lg(\bar x)}|\cM_n|$ we have
\end{enumerate}

\begin{enumerate}
\item[$(\ast\ast)$]
$\cM_n\models\varphi(\bar a)$ if and only if $\cM_n
\upharpoonright\cl^{k_{\varphi},m_{\varphi}}
(\bar a,\cM_n)\models\psi_\varphi(\bar a)$.
\end{enumerate}

\noindent 2) Moreover, if for simplicity we will consider
``$y\in\cl^{k,m}(\bar x,M)$'' as an atomic formula when computing
the
\footnote{  q.d stand for quantifier depth} q.d. of $\psi_{\varphi}$ {\em
then} we can demand: the number of alternation of quantifiers of $\psi
_\varphi$ is $\leq$ those of $\varphi$, more fully if $\varphi$ is a
$\Pi_n$ (or $\Sigma_n$) then so is $\psi_\varphi$.

\end{mainlemma}

\noindent \Proof  We shall ignore $(2)$, (which is not used and
is obvious if we understand the proof below). We prove the statement in part
($1$) by induction on $r=q.d.(\varphi(\bar x))$ and first note (by clause
(e) of Def.2.2 as ``$y \in cl^{k,m}(\bar x)$'' is f.o. definable in $\cK$)
that $(\ast)_{\varphi}$ implies 
\begin{enumerate}
\item[$(\ast)^+_{\varphi}$] in $(\ast)_{\varphi}$, possibly changing
$\psi_{\varphi}$ one can replace
$\cM_n\upharpoonright\cl^{{k_\varphi},m_{\varphi}}(\bar a,\cM_n)$
by any $N$ with
$\cl^{k_{\varphi},m_{\varphi}}(\bar a,\cM_n) \subseteq N \subseteq \cM_n$.
\end{enumerate}

\noindent Case 1 $\varphi$ atomic. Trivial [Proof of Case 1: If
$\varphi(\bar x)$ is an atomic formula
we let $m_{\varphi}=0$, $k_{\varphi}=0$ or whatever. So $\cl^{k_
\varphi,m_\varphi}(\bar a,\cM_n)=\bar a$ for our $k_\varphi$, $m_\varphi$.
Assume $\cM_n \models \varphi(\bar a)$ and we let $\psi_\varphi=\varphi$.
Now as $\bar a\subseteq \cM_n \restriction \cl^{k_\varphi,m_\varphi}(\bar a,
\cM_n)\subseteq \cM_n$ we have $\cM_n \models \varphi(\bar a)$ iff $cl^{k_
\varphi,m_\varphi}(\bar a,\cM_n) \models \psi_\varphi(\bar a)$ as required].\\

\noindent Case 2: $\varphi$ is a Boolean combination of atomic formulas and
the formulas of the form $\exists x \varphi(x,\bar y)$ with
q.d.$(\varphi')<r$. Clearly follows by case 3 and case 1. 

\noindent Case 3: $r>0$ and $\varphi (\bar x)=(\exists
y)\varphi_1 (\bar x,y)$. Let
$$
\begin{array}{ll}
(\ast)_1 \quad & m^\ast=m^\ast(k_{\varphi_1},m_{\varphi_1}, \lg(\bar x),
\ell'),\ \ k_\varphi=k^*(k_{\varphi_1},m_{\varphi_1},\lg(\bar x),\ell'),\\
& t=t(k_{\varphi_1},m_{\varphi_1},\lg(\bar x),\ell')
\end{array}
$$
be as guaranteed in Def.\ref{2.5} with $\ell'$ suitable (see its use below)
and let $m_{\varphi}:= m^\ast+m_{\varphi_1}$. Let $\psi^1_{\varphi_1}$ be
such that it witness $(\ast)_{\varphi_1}$, and let $\psi^2_{\varphi_1}$ be
such that it witness $(\ast)_{\varphi_1}$.

So it is enough to prove the following two statements:\\
\noindent {\em Statement 1:}\ \ \ There is $\psi^{1}_\varphi(\bar x)$ (f.o)
such that:
\begin{enumerate}
\item[$(\boxtimes)_1$] for every random enough $\cM_n$, for every $\bar a\in
{}^{\lg (\bar x)} |\cM_n|$ we have $(\alpha)_1 \Leftrightarrow (\beta)_1$
where:
\begin{enumerate}
\item[$(\alpha)_1$] $\cM_n\upharpoonright\cl^{k_\varphi,m_\varphi} (\bar a,
\cM_n)\models\psi^1_\varphi(\bar a)$
\item[$(\beta)_1$]  $\cM_n\models$ ``there is $b\in\cl^{k_\varphi,m^\ast}(\bar
a,\cM_n)$ such that $\varphi_1(\bar a,b)$ holds'' (i.e. b belongs to
a small enough neighbourhood of $\bar a)$.
\end{enumerate}
\end{enumerate}
\medskip

\noindent {\em Statement 2:}\ \ \ There is $\psi^2_\varphi(\bar x)$ (f.o) such
that:
\begin{enumerate}
\item[$(\boxtimes)_2$] for every random enough $\cM_n$ and for every 
$\bar a\in {}^{\lg(\bar x)} |\cM_n|$ we have $(\alpha)_2\Leftrightarrow
(\beta)_2$ where: 
\begin{enumerate}
\item[$(\alpha)_2$] $\cM_n\upharpoonright \cl^{k,m^\ast}(\bar a,\cM_n)\models
\psi^2_\varphi(\bar a)$
\item[$(\beta)_2$]  $\cM_n\models$ ``there is $b\in\cM_n\setminus
\cl^{k_\varphi,m^\ast}(\bar a,\cM_n)$ such that $\varphi_1(\bar a,b)$ holds''
(i.e. $b$ is far to $\bar a$)
\end{enumerate}
\end{enumerate}
(note: $(\beta)_1$, $(\beta)_2$ are complementary, but it is enough that
always at least one of them holds).

\noindent Note that as ``$y\in\cl^{k_\varphi,m^\ast}(\bar x)$'' is f.o
definable and $m{_{\varphi}} = m^{\ast}+m_{\varphi_1} \geq m^\ast$,
by \ref{2.1} and
clause (e), we can in $(\alpha)_2$ replace $m^\ast$ by
$m_\varphi$, changing $\psi^2_\varphi$ to $\psi^{2.5}_\varphi$.

\noindent Clearly these two statements are enough and $\psi^1_\varphi(\bar x)
\vee \psi^{2.5}_\varphi(\bar x)$ is as required.
\medskip

\noindent {\em Proof of statement 1:}

\noindent Easy, recalling that $k^\ast\geq k_{\varphi_1}$ by clause ($A$)
of Def.\ref{2.5}, by the induction hypothesis as
$$
\cl^{k_{\varphi_1},m_{\varphi_1}}(\bar ab,M_n)\subseteq\cl^{k_\varphi, m^\ast
+ m_{\varphi_1}}(\bar a,\cM_n)=\cl^{k_\varphi, m_\varphi}(\bar a,\cM_n)
$$
and by the fact that the closure is sufficiently definable.
\medskip

\noindent {\em Proof of statement 2:}

\noindent We will use a series of equivalent statements $\otimes_\ell$.
\begin{enumerate}
\item[$\otimes_1$] is $(\beta)_2$
\item[$\otimes_2$] there are $m_2\in [m_{\varphi_1},m^*]$, $m_1\le
m^\ast-m_2$,
 $b$, $B$ and $B^*$, $B'$ such that:
\end{enumerate}
\begin{enumerate}
\item[$(\alpha)$] $b\in\cM_n$, $b\notin\cl^{k_\varphi,m^\ast}(\bar a,\cM_n)$,
$\bar a
\subseteq B\subseteq\cl^{k_{\varphi_1},m_1}(\bar a,\cM_n)$, $|B|\le t$,
\item[$(\beta)$] $B^*=B\cup [\cl^{k_{\varphi_1},m_{\varphi_1}}(\bar
ab,\cM_n)\setminus
\cl^{k_{\varphi_1},m_2}(B,\cM_n)]$
[hence $B=B^*\cap\cl^{k_{\varphi_1},m_2}(B,\cM_n)$] and
\item[$(\gamma)$] $B\leq_s B' \in \cK_\infty$ and $B'=B^*$ or at least
$(B^*,b,c)_{c \in B}\equiv_{\ell'} (B',b,c)_{c \in B}$ (see \ref{2.5A}(4))
and
\item[$(\delta)$]
$$
\nonforkin{B^*}{\cl^{k_{\varphi_1}, m_2}(B,\cM_n)}_{B}^{\cM_n}
$$
\item[$(\vep)$] $(B',B,\bar ab,B^* \cap \cl^{k,m_0}(\bar ab,\cM_n),
k,m_0,m_2)$
is good,
\item[$(\zeta)$] for $m \leq m_{\varphi_1}$ we have $\cl^{k_{\varphi_1}, m}
(\bar ab, B^*)= B^* \cap \cl^{k_{\varphi_1}, m}(\bar ab,\cM_n)$
\end{enumerate}
and
\begin{enumerate}
\item[$\oplus_2$] $\cM_n\models\varphi_1(\bar a,b)$
\end{enumerate}

\begin{enumerate}
\item[$(\ast)_2$] $\otimes_1 \Leftrightarrow \otimes_2$
\end{enumerate}
Why?  The implication $\Leftarrow$ is trivial
as $\oplus_2$ is included in $\otimes_2$, the implication $\Rightarrow$
holds by clause (A) in the definition of almost nice \ref{2.5}, except
$b \notin cl^{k_{\varphi},m^{\ast}}(\bar a,\cM_n)$ which is explicitly
demanded in $(\beta)_2$.

\begin{enumerate}
\item[$\otimes_3$] like $\otimes_2$ but replacing $\oplus_2$ by
\item[$\oplus_3$] $\cM_n\upharpoonright\cl^{k_{\varphi_1},m_{\varphi_1}}(\bar
ab,\cM_n)\models\psi^1_{\varphi_1}(\bar a,b)$.
\end{enumerate}

\begin{enumerate}
\item[$(\ast)_3$]  $\otimes_2 \Leftrightarrow \otimes_3$
\end{enumerate}
Why? By the induction hypothesis.

\begin{enumerate}
\item[$\otimes_4$] like $\otimes_3$ replacing $\oplus_3$ by
\begin{enumerate}
\item[$\oplus_4$] $\cM_n\upharpoonright [B^* \cup \cl^{k_{\varphi_1},m_2}
(B,\cM_n)]\models\psi^2_{\varphi_1}
(\bar a,b)$.
\end{enumerate}
\end{enumerate}

\begin{enumerate}
\item[$(\ast)_4$] $\otimes_3\Leftrightarrow\otimes_4$
\end{enumerate}
Why? By $(\ast)^+_{\varphi_1}$ in the beginning of the proof, the
definition of $B^*$ and
the choice of $\psi^2_{\varphi_1}$ (Let $\otimes_3$ be true. As by the
choice of $B^*$, $B$ above,
$\cl^{k_{\varphi_1},m_{\varphi_1}}(\bar ab,\cM_n) \cup \cl^{k_{\varphi_1},
m_2}(B,\cM_n)= B^* \cup \cl^{k_{\varphi_1},m_2}(B,\cM_n) \subseteq
\cM_n$ we have $\cM_n \models \varphi_1(\bar a,b)$ iff\\
 $B^* \cup \cl^{k_{\varphi_1},m_2}(B,\cM_n) \models \psi^2_{\varphi_1}
(\bar ab)$ by $(*)^+_{\varphi_1})$.
So $(\ast)_4$ holds.

For notational simplicity we assume $B\neq \emptyset$, and similarly
assume
$\bar a$ is with no repetition and we shall apply the lemma \ref{2.6A}
several times.

First we apply \ref{2.6A} to the case $s=t$, $\bar z=\langle z_1,\ldots,
z_t\rangle$, $\bar z^1=\langle z^1_1, z^1_2\rangle$, $\bar z^2$ empty and
the formula ``$z^1_2\in\cl^{k_{\varphi_1},m_{\varphi_1}}(\bar z,z^1_1)$''
and get $i^*_{1,m}\in\bbn$ and formulas $\theta^1_{1,m,i}(\bar z,z^1_1,
z^1_2)$ and $\theta^2_{1,m,i}(\bar z)$ for $i<i^*_{1,m}$. Let 
$$
u^*_1=\{(m,i): m\leq m_{\varphi_1},i<i^*_{1,m}\}.
$$
Second for $m\leq m_{\varphi_1}$ we apply \ref{2.6A} to the case $s=t$,
$\bar z^2
=\langle z^2_1\rangle$, $\bar z^1=\langle z^1_1\rangle$, $\bar z=\langle
z_1,\ldots,z_t\rangle$ and the formula ``$z^2_1\in\cl^{k_{\varphi_1},
m}(\bar z,z^1_1)$'' and get $i^*_{2,m}\in\bbn$ and formulas
$\theta^1_{2,m,i}(\bar z,\bar z^1)$ and $\theta^2_{2,m,i}(\bar z,z^2_1)$, for
$i<i^*_{2,m}$.

Let $\tau'=\tau_{\cK} \cup \{P_1,P_2\}$, with $P_1$, $P_2$
new unary predicates: for $\theta\in\cL [\tau'_{\cK}]$ let
$\theta^{[P_\ell]}$ be $\theta$ restricting the quantifiers to $P_\ell$.
Let $\psi^*=\psi^*_1\wedge\psi^*_2\wedge\psi^*_3$ where
$$
\begin{array}{ll}
\psi^*_1 =: &\psi^2_{\varphi_1}(z_1,\ldots,z_{\lg(\bar x)},z^1_1)\\
\psi^*_2 =: &\bigwedge\limits_{m\leq m_{\varphi_1}} (\forall y)\Bigl[y\in
\cl^{k_{\varphi_1}, m_{\varphi_1}}(\{z_1,\ldots,z_{\lg(\bar x)}, z^1_1 \})\\
\ &\qquad\qquad\qquad\equiv \bigl(\psi^{*,1}_{2,m}(z_1,\ldots,z_t,z^1_1,y)\vee
\psi^{*,2}_{2,m}(z_1,\ldots,z_t,z^1_1,y)\bigr)\Bigr],
\end{array}
$$
where
$$
\psi^{*,1}_{2,m}(z_1,\ldots,z_t,z^1_1,y)=:\bigvee_{i<i^*_{1,m}}(\theta^1_{1,
m, i}(z_1,\ldots,z_t,z^1_1,y)^{[P_1]}\wedge \theta^2_{1, m, i}(z_1, \ldots,
z_t)^{[P_2]})
$$
$$
\psi^{*,2}_{2,m}(z_1,\ldots,z_t,z^1_1,y)=:
\bigvee_{i<i^*_{2, m}}(\theta^2_{1,
m,i}(z_1,\ldots,z_t,z^1_1)^{[P_1]}\wedge\theta^2_{2,m,i}(z_1,\ldots,z_t,
y)^{[P_2]})
$$
and let

\begin{multline*}\psi^*_3=:(\forall y)\Bigl(P_1(y)\equiv[\bigvee^t_{\ell=1}
y=z_{\ell}\vee\\
y\in\cl^{k_{\varphi_1},m_{\varphi_1}}(\{z_1,\ldots,z_{\lg(\bar x)},z^1_1\}\\
\vee y\notin \cl^{k_{\varphi_1,m_2}} \{ z_1,...,z_
{lg(\bar x)} \}] \Bigr).
\end{multline*}

So we have defined $\psi^*$. Now we apply \ref{2.6A} the third time, with the
vocabulary $\tau_{\cK}\cup\{P_1,P_2\}$ to the case $s=t$, $\bar z^2$ empty,
$\bar z^1=\langle z^1_1\rangle$, $\bar z=\langle z_1,\ldots,z_\ell\rangle$,
and $\psi(\bar z,\bar z^1,\bar z^2)=\psi(\bar z,z^1_1)=\psi^*(\langle z_1,
\ldots,z_{\lg\bar (x)}\rangle,z^1_1)$ and get $i^*$, $\theta^1_{3,i}(\bar z,
\bar z^1)$ and $\theta^2_{3,i}(\bar z,\bar z^2)$ as there.
Let

\begin{enumerate}
\item[$\otimes_5$] like $\otimes_4$ but replacing $\oplus_4$ by
\begin{enumerate}
\item[$\oplus_5$] letting $c_1,\ldots,c_t$ list $B$ possibly with repetitions
but such that $\langle c_1,\ldots,c_{\lg(\bar x)}\rangle=\bar a$ and letting
$$
P^*_1=B^*\mbox{ and }P^*_2=\cl^{k_{\varphi_1},m_2}(\{c_1,\ldots,c_t\},\cM_n)
$$
we have
\begin{enumerate}
\item[$(\ast)$] $(\cM_n\restriction (P^*_1\cup P^*_2),P^*_1,P^*_2)\models
\psi^*[c_1,\dots,c_t,b]$  (the model is a $\tau'$-model).
\end{enumerate}
\end{enumerate}
\end{enumerate}
Now
\begin{enumerate}
\item[$(*)_5$] $\otimes_4 \Leftrightarrow \otimes_5$.
\end{enumerate}
Why? Look at what the statements mean recalling $\nonforkin{\cM_n
\restriction P^*_1}{\cM_n \restriction P^*_2}_{B}^{\cM_n}$. Next let
\begin{enumerate}
\item[$\otimes_6$] like $\otimes_5$ but replacing $\oplus_5$ by
\begin{enumerate}
\item[$\oplus_6$] letting $c_1,\ldots,c_t$ list $B$ possibly with repetitions
but such that $\langle c_1,\ldots,c_{\lg(\bar x)}\rangle=\bar a$ and letting
$$
P^*_1= B^*\mbox{ and }P^*_2=\cl^{k_{\varphi_1},m_1}(\{c_1,\ldots,c_t\},\cM_n)
$$
{\em there is} $i<i^*$ such that:
\begin{enumerate}
\item[(i)]  $(\cM_n\restriction P^*_1,P^*_1,P^*_2\cap P^*_1)\models
\theta^1_{3, i}[\langle c_1,\ldots,c_t\rangle,b]$,
\item[(ii)] $(\cM_n\restriction P^*_2, P^*_1\cap P^*_2, P^*_2)\models
\theta^2_{3,i}[\langle c_1,\ldots,c_t\rangle]$.
\end{enumerate}
\end{enumerate}
\end{enumerate}
Now
\begin{enumerate}
\item[$(*)_6$] $\otimes_5 \Leftrightarrow \otimes_6$
\end{enumerate}
Why? By the choice of $\theta^1_{3,i}$, $\theta^2_{3,i}$ ($i<i^*$).

However in the two $\tau'$-models appearing in $\oplus_6$, the predicates
$P_1$, $P_2$ are interpreted in a trivial way: as the whole universe of the
model or as $\{c_1,\ldots,c_t\}$.

 So let:
\begin{enumerate}
\item[(a)]  $\theta^1_{4,i}(z_1,\ldots, z_t, y)$ be $\theta^1_{3,i}(z_1,
\ldots,z_t,y)$ with each atomic formula of the form $P_1(\sigma)$ or
$P_2(\sigma)$ being replaced by $\sigma=\sigma$ or $\bigvee^{t}_{r=1}
\sigma=z_r$ respectively,
\item[(b)] $\theta^2_{4,i}(z_1,\ldots,z_t)$ be $\theta^2_{3,i}(z_1,\ldots,
z_t)$ with each atomic formula of the form $P_1(\sigma)$ or $P_2(\sigma)$
being replaced by $\bigvee^t_{r=1}\sigma=z_r$ or $\sigma=\sigma$ respectively.
\end{enumerate}

So let (recall $B'$ is mentioned in $\otimes_2$, a ``replacement'' to $B^*)$
\begin{enumerate}
\item[$\otimes_7$] like $\otimes_6$ but replacing $\oplus_6$ by
\begin{enumerate}
\item[$\oplus_7$] letting $c_1,\ldots,c_t$ list $B$ possibly with repetitions
but such that $\langle c_1,\ldots,c_{\lg\bar x}\rangle=\bar a$, {\em there is}
$i<i^*$ such that
\begin{enumerate}
\item[(i)]  $\cM_n\restriction B'\models\theta^1_{4,i}[\langle c_1,\ldots,
c_t\rangle,b]$ and
\item[(ii)] $\cM_n\restriction\cl^{k_{\varphi_1},m_2}(\langle c_1,\ldots,
c_t\rangle,\cM_n) \models\theta^2_{4,i}(\langle c_1,\ldots,c_t\rangle)$.
\end{enumerate}
\end{enumerate}
\end{enumerate}

\begin{enumerate}
\item[$(\ast)_7$] $\otimes_6\Leftrightarrow \otimes_7$
\end{enumerate}
Why? By the choice of the $\theta^1_{4,i}$, $\theta^2_{4,i}$ and the property
of $B'$ (stated in $\otimes_2$).

Let $\cP=\{(N,c_1,\ldots,c_t): N\in\cK_\infty$, with the set of elements
 $\{c_1,\ldots,c_t\}\}$. Let $\{(N_j,c_1^j,\ldots,c^j_t):j<j^*\}$ list the
 members of $\cP$ up to isomorphism, so with no two isomorphic. For every
 $j<j^*$ and $i<i^*$ choose if possible $(N_{j,i},c^j_1,\ldots,c^j_t,b^j_i)$
 such that: 
\begin{enumerate}
\item[(i)]   $N_j\leq_s N_{j,i}$ (in $\cK_\infty$),
\item[(ii)]  $b^j_i\in N_{j,i}\setminus N_j$,
\item[(iii)] $N_{j,i}\models\theta^1_{4,i}(\langle c^j_1,\ldots,c^j_t\rangle,
b^j_i)$  and
\item[(iv)]  $(N_{j,i},B,\{c^j_i:i=1,\ldots,\lg(\bar
x)\}\cup\{b^j_i\},k,m_0,m_2)$ is good for $\GK$.
\end{enumerate}
Let
$$
w=\{(i,j):i<i^*,j<j^*\ \mbox{ and }\ (N_{j,i},c^j_1,\ldots,c^j_t,b^j_i)\mbox{
is well defined}\}.
$$
Let
\begin{enumerate}
\item[$\otimes_8$] there are $m_2\leq m^*$, $m_1\leq m^* - m_2$, such that
$m_2\geq m_{\varphi_1}$ and, there are $b$, $B$ such that:

$\bar a\subseteq B\subseteq\cl^{k^*,m_2}(\bar a,\cM_n)$, $|B|\leq
t(k_{\varphi_1},m_{\varphi_1},\lg(\bar x))$, $b\notin\cl^{k^*,m^*}(\bar
a,\cM_n)$, $b\in\cM_n$, and
\begin{enumerate}
\item[$\oplus_8$] for some $c_1,\ldots,c_t$ listing $B$ such that $\bar
a=\langle c_1,\ldots,c_{\lg\bar x}\rangle$ {\em there are} $i<i^*$, $j<j^*$
such that $(i,j) \in w$ and:
\begin{enumerate}
\item[(i)] $(\cM_n\restriction B,c_1,\ldots,c_t)\cong(N_j,c^j_1,\ldots,c^j_t)$
i.e. the mapping $c^j_1 \mapsto c_1$, $c^j_2 \mapsto c_2$ embed $N_j$ into
$\cM_n$,
\item[(ii)] $\cM_n\restriction\cl^{k_{\varphi_1},m_2}(B,\cM_n)\models
\theta^2_{4,i}(\langle c_1,\ldots,c_t\rangle)$
\end{enumerate}
\end{enumerate}
\end{enumerate}

\begin{enumerate}
\item[$(*)_8$] $\otimes_7 \Leftrightarrow \otimes_8$
\end{enumerate}
Why? For proving $\otimes_7 \Rightarrow \otimes_8$ let $c_1,\ldots,c_t$ as
well as $i<i^*$ be as in $\oplus_7$, let $j<j^*$ be such that $(\cM_n
\restriction B,c_1,\ldots,c_t)\cong(N_j,c^j_1,\ldots,c^j_t)$. The main point
is that $B'$ exemplifies that $(i,j)\in w$.\\
For proving $\otimes_8 \Rightarrow \otimes_7$ use the definition of goodness
in clause $(\vep)$ (see $\otimes_2$ and Def. in \ref{2.7new}(4).

We now have finished as $\otimes_8$ can be expressed as a f.o formula
straightforwardly. So we have carried the induction hypothesis on the
quantifier depth thus finishing the proof.  \hfill\rqed$_{\ref{2.6}}$

\begin{lemma}
\label{2.14}
1) Assume $(\GK,\cl)$ is simply almost nice
and $\cl $ is f.o. definable. Let $\varphi(\bar x)$ be a f.o.
formula in the vocabulary $ \tau_{\mathcal K} $. {\em Then} for some
$k=k_\varphi$ and f.o. formula $\psi_\varphi(\bar x)$ we have:
\begin{enumerate}
\item[$(\ast)$] for every random enough $\cM_n$ and $\bar a\in
{}^{\lg(\bar x)}
|\cM_n|$
\item[$(\ast\ast)$] $\quad\cM_n\models\varphi(\bar a)$ if and only if $\cM_n
\upharpoonright\cl^{k_\varphi}(\bar a,\cM_n)\models\psi_\varphi(\bar a)$
\end{enumerate}
\noindent 2) The number of alternation of quantifiers of $\psi_\varphi$ in
$(1)$ is $\leq$ the number of alternation of quantifiers of $\varphi$ if
we consider ``$y\in\cl^{k,m}(\bar x,M)$'' as atomic. More fully, if $\varphi$
is $\Pi_n$ (or $\Sigma_n)$ then $\psi_\varphi$ is.
\end{lemma}

\begin{remark}
\label{2.20new}
\begin{enumerate}
\item Of course we do not need to assume that closure operation is definable,
it is enough if there is a variant  $\cl_*$  which is definable and for
every $k,m$ there are $ k^1, m^1, k_2,m_2 $ such that always $ \cl^{k,m}
(A,M)   \subseteq \cl_*^{k^1,m^1} (A,M ) \subseteq \cl^{k^2,m^2}(A,M) $.
\item Similarly in \ref{2.6} (using Def.\ref{2.7A}).
\item We can weaken ``simply almost nice'' as in Def.\ref{2.12Anew}(1)
and still part $(1)$ is true, with essentially the same proof.
\item The proof of \ref{2.14} is somewhat simpler than the proof of
\ref{2.6}. 
\end{enumerate}
\end{remark}

\Proof 1) We prove the statement by induction on $r=q.d.(\varphi(\bar x))$.
First note (by clause (e) of \ref{2.1})
\begin{enumerate}
\item[$(\ast)^+_{\varphi}$] in  $(\ast)$ (of \ref{2.14}, possibly changing
$\psi_\varphi$) one can
replace $\cM_n\upharpoonright\cl^{k_\varphi}(\bar a,\cM_n)$ by any $N$ with
$\cl^{k_\varphi}(\bar a,\cM_n)\subseteq N\subseteq\cM_n$.
\end{enumerate}

\noindent Case 1: Let $\varphi$ be atomic.\ \ \ \ \ \ Trivial.
\medskip

\noindent Case 2: $\varphi$ a Boolean combination of atomic formulas and
formulas of quantifier depth $<r$. Clearly follows by case 3 and case 1. 
Trivial.
\medskip

\noindent Case 3: $r>0$ and $\varphi(\bar x)=(\exists y)\varphi_1
(\bar x,y)$. Let (the functions are from \ref{2.12new}(1))
$$
m^\ast=m^\ast(k_{\varphi_1},\lg (\bar x),\ell'),\quad k^*=k^*(k_{\varphi_1},
\lg(\bar x),\ell'),\quad t=t(k_{\varphi_1},\lg(\bar x),\ell')
$$
with $\ell'$ suitable (just the quantifier depth of $\psi^2_{\varphi_1}$
defined below) and let $k_\varphi$ be\footnote{if we change clause (A)
of \ref{2.12new}(1)
a little, $k_\varphi=k^*$ will be O.K.: instead of assuming
$b\notin\cl^{k^*,m^*}(\bar a,\cM_n)$ assume just $\cl^k(\bar ab,\cM_n)
\not\subseteq
\cl^{k^*,m^*}(\bar a,\cM_n)$. Allowing to increase $m^*$, the two versions are
equivalent. $m^{**}=m^{**}(k,\ell,\ell')=m^*(k,\ell,\ell')+k$. Now by $2.4(3)$
we have $b \in \cl^{k^*,m^*}(\bar a, \cM_n)$ and $c \in \cl^k(\bar ab,\cM_n)
\Rightarrow c \in \cl^{m^*+k}(\bar a,\cM_n)=\cl^{m^{**}}(\bar a,\cM_n)$ hence
$b \in \cl^{k^*,m^*}(\bar a,\cM_n) \Rightarrow \cl^k(\bar ab,\cM_n)
\subseteq \cl^{k,m^{**}}(\bar a,\cM_n)$ hence $\cl^k(\bar ab,\cM_n)
\nsubseteq \cl^{k,m^{**}}(\bar a,\cM_n) \Rightarrow b \notin \cl^{k^*,m^*}
(\bar a,\cM_n)$, so our new assumption for $m^{**}$ implies  our old for
$m^*$. Of course our new assumption for $m^*$ implies  the old for $m^*$.
See section 3 where this is done.
}   such that:
$$
(\ast)_1\quad
|A|\leq\lg(\bar x)+1\ \&\ A\subseteq N\in {\mathcal K}\ \Rightarrow\
\cl^{k_{\varphi_1}}(\cl^{k^*,m^*}(A, N),N)\subseteq\cl^{k_\varphi}(A,N).
$$
Let $\psi^1_{\varphi_1}(\bar{x},y)$ be such that it witness
$(\ast)_{\varphi_1}$ holds, and let $\psi^2_{\varphi_1}(\bar{x},y)$ be
such that it witness $(\ast)^+_{\varphi_1}$.

It is enough to prove the following two statements (see below):
\medskip

\noindent {\em Statement 1:}  There is $\psi^{1}_\varphi (\bar x)$ (f.o.) such
that:
\begin{enumerate}
\item[$(\boxtimes)_1$] for every random enough $\cM_n$, for every $\bar a\in
{}^{\lg(\bar x)}|\cM_n|$ we have $(\alpha)_1\Leftrightarrow (\beta)_1$ where:
\begin{enumerate}
\item[$(\alpha)_1$] $\cM_n\upharpoonright\cl^{k_{\varphi}}(\bar a,\cM_n)
\models \psi^1_\varphi(\bar a)$
\item[$(\beta)_1$] $\cM_n\models$ ``there is $b\in\cl^{k^*,m^\ast}(\bar
a,\cM_n)$ such that $\varphi_1(\bar a,b)$ holds.''
\end{enumerate}
\end{enumerate}
\medskip

\noindent {\em Statement 2:}  There is $\psi^2_\varphi(\bar x)$ (f.o) such
that:
\begin{enumerate}
\item[$(\boxtimes)_2$] for every random enough $\cM_n$ and for every $\bar
a\in {}^{\lg(\bar x)} |\cM_n|$ we have $(\alpha)_2\Leftrightarrow (\beta)_2$
where: 
\begin{enumerate}
\item[$(\alpha)_2$] $\cM_n\upharpoonright\cl^{k^*,m^\ast}(\bar a,\cM_n)
\models\psi^2_\varphi(\bar a)$
\item[$(\beta)_2$] $\cM_n\models$ ``there is $b\in\cM_n\setminus
\cl^{k^*,m^\ast}(\bar a,\cM_n)$ such that $\varphi_1 (\bar a,b)$ holds''
\end{enumerate}
\end{enumerate}
(note: $(\beta)_1$, $(\beta)_2$ are complementary, but it is enough that
always at least one holds).

\noindent Note that as ``$y\in\cl^{k^*,m^\ast}(\bar x)$'' is f.o. definable,
by \ref{2.1}, clause (e) and the choice of $k_\varphi$ we can in $(\alpha)_2$
replace $\cl^{k^*,m^\ast}$ by $\cl^{k_\varphi}$, changing $\psi^2_\varphi$ to
$\psi^{2.5}_\varphi$; (just as from$(*)_\psi$ we have deduce
$(*)^+_{\varphi}$.

\noindent Clearly these two statements are enough as if
$\psi^{1.5}_\varphi$ express $(\alpha)_1$
then $\psi^{1.5}_\varphi
(\bar x) \vee \psi^{2.5}_\varphi (\bar x)$ is as required.
\medskip

\noindent {\em Proof of statement 1:}

\noindent Easily, by the induction hypothesis as
$$
\cl^{k_{\varphi_1}}(\bar ab, M_n)\subseteq\cl^{k_{\varphi_1}}(\cl^{k^*,
m^\ast}(\bar a,\cM_n),\cM_n)\subseteq\cl^{k_\varphi}(\bar a,\cM_n)
$$
and by the fact that the closure is sufficiently definable. So in this case
$\psi_\varphi(\bar{a})$ can be chosen as $(\exists y)\psi^2_{\varphi_1}
(\bar{a},y)$.
\medskip

\noindent {\em Proof of statement 2:}

\noindent We will use a series of equivalent statements $\otimes_\ell$.
\begin{enumerate}
\item[$\otimes_1$] is $(\beta)_2$
\item[$\otimes_2$] there are $b$, $B$ and
$B^*$, $B'$ such that:

\item[$(\alpha)$] $b\in\cM_n$, $b\notin\cl^{k^*,m^\ast}(\bar a,\cM_n)$,

\item[$(\beta)$] $\bar a\subseteq B\subseteq\cl^{k^*,m^*}(\bar a,\cM_n)$,
moreover $\cl^{k_{\varphi_1}}(B,\cM_n) \subseteq\cl^{k^*,m^*}
(\bar a,\cM_n)$, and $|B| \leq t$,

\item[$(\gamma)$] $B^*\supseteq B\cup [\cl^{k_{\varphi_1}}(\bar ab,\cM_n)
\setminus\cl^{k_{\varphi_1}}(B,\cM_n)]$ and

\item[$(\delta)$] $B\leq_s B'\in\cK_\infty$ and: $B'=B$ or just
$(B^*,b,c)_{c\in B}\equiv_{\ell'}
(B',b,c)_{c\in B}$ (see \ref{2.5A}(4)) and

\item[$(\vep)$] $\nonforkin{B^*}{\cl^{k_{\varphi_1}}(B, \cM_n)}_{B}^{\cM_n}$
(and so $B=B^*\cap\cl^{k_{\varphi_1}}(B,\cM_n)$) and

\item[$(\zeta)$] $(B', B, \bar ab, k)$ is simply good

\item[$(\eta)$] $\cl^k(\bar ab,B^*) \setminus B=B^* \cap \cl^k(\bar ab,
\cM_n) \setminus  \cl^k(B, {\mathcal M}_n ) $,
actually this follows from clauses $(\epsilon), (\beta)$, and
\item[$\oplus_2$] $\cM_n \models \varphi_1(\bar a,b)$
\end{enumerate}

\begin{enumerate}
\item[$(\ast)_2$] $\otimes_1 \Leftrightarrow \otimes_2$
\end{enumerate}
Why?  The implication $\Leftarrow$ is trivial
as $\oplus_2$ is included in $\otimes_2$, the implication $\Rightarrow$
holds by clause (A) in the definition \ref{2.12new} of simply almost nice.

\begin{enumerate}
\item[$\otimes_3$] like $\otimes_2$ but replacing $\oplus_2$ by
\begin{enumerate}
\item[$\oplus_3$] $\cM_n\upharpoonright\cl^{k_{\varphi_1}}(\bar ab,M_n)
\models\psi^1_{\varphi_1}(\bar a,b)$.
\end{enumerate}
\end{enumerate}

\begin{enumerate}
\item[$(\ast)_3$]  $\otimes_2 \Leftrightarrow \otimes_3$
\end{enumerate}
Why?  By the induction hypothesis and our choices.   

\begin{enumerate}
\item[$\otimes_4$] like $\otimes_3$ replacing $\oplus_3$ by
\begin{enumerate}
\item[$\oplus_4$] $\cM_n \upharpoonright[B^*
\cup\cl^{k_{\varphi_1}}(B,\cM_n)]\models\psi^2_{\varphi_1}(\bar a,b)$.
\end{enumerate}
\end{enumerate}

\begin{enumerate}
\item[$(\ast)_4$] $\otimes_3 \Leftrightarrow \otimes_4$
\end{enumerate}
Why? By $(\ast)^{+}_{\varphi_1}$ in the beginning of the proof, the
requirements on $B^*$ and the choice of $\psi^2_{\varphi_1}$.

For notational simplicity we assume $B\neq \emptyset$, and similarly assume
$\bar a$ has no repetitions and apply the lemma \ref{2.6A} with the vocabulary
$\tau_\GK$ to the case $s=t$, $\bar z^2$ empty, $\bar z^1=\langle
z^1_1
\rangle$, $\bar z=\langle z_1,\ldots,z_t\rangle$, and $\psi(\bar z,\bar
z^1, \bar z^2)=\psi(\bar z,z^1_1)=\psi_{\varphi_1}^2(\langle z_1,\ldots,
z_{\lg\bar x}\rangle,z^1_1)$ and get $i^*$, $\theta^1_{i}(\bar z,\bar z^1)$
and $\theta^2_{i}(\bar z)$ for $i<i^*$ as there; in particular the quantifier
depth of $\theta^1_i$, $\theta^2_i$ for $i<i^*$ is at most the quantifier
depth of $\psi^2_{\varphi_1}$.

\noindent Next let
\begin{enumerate}
\item[$\otimes_5$] like $\otimes_4$ but replacing $\oplus_4$ by
\begin{enumerate}
\item[$\oplus_5$] letting $c_1,\ldots,c_t$ list $B$ possibly with repetitions
but such that $\langle c_1, \ldots, c_{\lg(\bar x)}\rangle=\bar a$,
$i<i^*$ such that:
\begin{enumerate}
\item[(i)]  $B^*\models\theta^1_{i}[\langle c_1,\ldots,c_t\rangle,b]$
\item[(ii)] $\cl^k(B,\cM_n)\models\theta^2_{i}[\langle c_1,\ldots,c_t\rangle]$
\end{enumerate}
\end{enumerate}
\end{enumerate}

Now
\begin{enumerate}
\item[$(*)_5$] $\otimes_4 \Leftrightarrow \otimes_5$
\end{enumerate}
Why? By the choice of $\theta^1_{i}$, $\theta^2_{i}$ for $i<i^*$, so
by lemma \ref{2.6A}.

Let $\cP=\{(N,c_1,\ldots,c_t): N\in\cK_\infty$, with the set of elements
$\{c_1,\ldots,c_t\}\}$. Let $\{(N_j,c_1^j,\ldots,c^j_t): j<j^*\}$ list the
members of $\cP$ up to isomorphism, so with no two isomorphic. For every
$j<j^*$ and $i<i^*$ choose if possible $(N_{j,i},c^j_1,\ldots,c^j_t,b^j_i)$
such that:
\begin{enumerate}
\item[(i)]   $N_j\leq_s N_{j,i}$ (in $\cK_\infty$),
\item[(ii)]  $b^j_i\in N_{j,i}\setminus N_j$,
\item[(iii)]  $N_{j,i}\models\theta^1_{i}(\langle c^j_1,\ldots,c^j_t\rangle,
b^j_i)$ and
\item[(iv)]   $(N_{j,i},\{c^j_1,\ldots,c^j_t\},\{c^j_1,\ldots,c^j_{\lg
\bar x}, b^j_i\},k)$ is simply good for $\GK$.
$$
w=\{(i,j): i<i^*,j<j^*\mbox{ and }(N_{j,i},c^j_1,\ldots,c^j_t,b^j_i)\mbox{
is well defined}\}.
$$
\end{enumerate}
Let

$\otimes_6$ like $\otimes_5$ replacing $\oplus_5$ by

$\oplus_6$ like $\oplus_5$ adding
\begin{enumerate}
\item[(iii)] for some $j$, $(i,j) \in w$ and $(B,c_1,...,c_t) \cong N_{j,i}$
\end{enumerate}

\begin{enumerate}
\item[$(\ast)_6$] $\otimes_5 \Leftrightarrow \otimes_6$
\end{enumerate}
Why? By the definition of $w$.

Let
\begin{enumerate}
\item[$\otimes_7$] there is $B$ such that: $b\in\cM_n$, $\bar a\subseteq
B\subseteq\cl^{k^*,m^*}(\bar a,\cM_n)$, $\cl^{k_{\varphi_1}}(B,\cM_n)
\subseteq\cl^{k^*,m^*}(\bar a,\cM_n)$, $|B|\leq t$, and
\begin{enumerate}
\item[$\oplus_7$] for some $c_1,\ldots,c_t$ listing $B$ such that $\bar a=
\langle c_1,\ldots,c_{\lg \bar x}\rangle$

{\em there are} $i<i^*$, $j<j^*$ such that $(i, j)\in w$ and:
\begin{enumerate}
\item[(i)] $(\cM_n\restriction B,c_1,\ldots,c_t)\cong (N_j,c^j_1,\ldots,
c^j_t)$ i.e. the mapping \\
$c^j_1\mapsto c_1$, $c^j_2\mapsto c_2$ embeds $N_j$
into $\cM_n$,
\item[(ii)] $\cM_n\restriction\cl^{k_{\varphi_1}}(B,\cM_n)\models\theta^2_{i}
(\langle c_1,\ldots,c_t\rangle)$
\end{enumerate}
\end{enumerate}
\end{enumerate}

\begin{enumerate}
\item[$(*)_7$] $\otimes_6 \Leftrightarrow \otimes_7$
\end{enumerate}
Why? For proving $\otimes_6 \Rightarrow\otimes_7$ let $c_1,\ldots,c_t$ as well
as $i<i^*$, $j<j^*$ be as in $\oplus_6$, let $j<j^*$ be such that
$(\cM_n\restriction
B,c_1,\ldots,c_t)\cong (N_j,c^j_1,\ldots,c^j_t)$. The main point is that $B'$
exemplifies that $(i,j)\in w$ (remember: $B'$ is from $\otimes_2$, and if
$B^*\in\cK_\infty$, we normally could have chosen $B'=B^*$).\\
For proving $\otimes_7 \Rightarrow \otimes_6$ use definition of simply
good tuples in Definition \ref{2.5B}(1).

We now have finished as $\otimes_7$ can be expressed as a f.o. formula
straightforwardly. So we have carried the induction hypothesis on the
quantifier depth thus finishing the proof.

\noindent 2) Similar
\hfill\rqed$_{\ref{2.14}}$

\begin{conclusion}
\label{2.7}
\begin{enumerate}
\item Assume $(\GK,\cl)$ is almost nice or simply almost nice and $\cl$ is
f.o. definable.

{\em Then}: $\GK$ satisfies the 0-1 law \ \ \ {\em iff}\ \ \ for any $k$, $m$
we have
\begin{enumerate}
\item[$(\ast)_{k,m}$] $\langle \cM_n \restriction \cl^{k,m}(\emptyset):
n<\omega\rangle\ \mbox{ satisfies the 0-1 law.}$
\end{enumerate}
\item Similarly with convergence and the very weak $0-1$ law.
\end{enumerate}
\end{conclusion}

\Proof 1) We first prove the ``only if''. There is a f.o. formula
$\theta(x)$ such that for every random enough $\cM_n$, $\theta(x)$ define
$\cl^{k,m}(\emptyset,\cM_n)$. Hence for every f.o. sentence $\varphi$ there
is a f.o. sentence $\psi_{\varphi}$ which is the relativization of $\varphi$
to $\theta(x)$, call it $\psi_\varphi$; hence, for every model $M \in \cK$,
$M \models \psi_\varphi\Leftrightarrow M \restriction \{a: M \models
\theta[a] \} \models \varphi$. Now  for every random enough $\cM_n$ we have
$a \in \cM_n \Rightarrow \cM_n \models \theta[a]\Leftrightarrow a \in
\cl^{k,m}(\emptyset,\cM_n)$, hence together
$$
\cM_n\models\psi_{\varphi}\Leftrightarrow \cM_n\restriction\cl^{k,m}
(\emptyset,\cM_n)
\models\varphi.
$$
As we are assuming that $\GK$ satisfies 0-1 law, for some truth value $\bt$
for every random enough
$\cM_n$
$$
\cM_n\models\mbox{``}\psi_{\varphi}\equiv\bt\mbox{''}
$$
hence (as required)
$$
\cM_n\restriction\cl^{k,m}(\emptyset,\cM_n)\models\mbox{``}\varphi=
\bt\mbox{''}.
$$
The other direction is similar by the main lemma \ref{2.6} when $(\cK,\cl)$
is almost nice, \ref{2.14} when $(\cK,\cl)$ is simply almost nice.

\noindent 2) Similar, so left to the reader. \hfill\rqed$_{\ref{2.7}}$

\begin{definition}
\label{2.20}
\begin{enumerate}
\item The tuple $(N, \bar b, \psi(\bar x), \langle B_0, B_1\rangle, k, k_1)$
is {\em simply$^*$ good} for $(\GK, \cl)$ if: $B_0$, $B_1\leq N\in
\cK_\infty$, $\cl^k(B_0, N)\subseteq B_1$, $\bar b\in {}^{\lg\bar
x}N$, $\psi(\bar x)$ a f.o. formula
and $k$, $k_1\in \bbn$ and for every random enough $\cM_n$, for every
$\bar b' \in {}^{\lg\bar x}(\cM_n)$ such that $\cM_n \restriction
\cl^{k_1}(\bar b', \cM_n) \models \psi(\bar b')$, letting $B'=\cM_n
\restriction \Rang (\bar b')$,
there is an embedding $g$ of $N$ into $\cM_n$ such that
\begin{enumerate}
\item[(i)] $g(\bar b)=\bar b'$
\item[(ii)] $g(N)\cap \cl^{k_1}(\bar b', \cM_n)= B'$
\item[(iii)] $g(N)\nonforkin{}{}_{B'}^{} \cl^{k_1}(\bar b', \cM_n)$
\item[(iv)] $\cl^{k}(g(B_0), \cM_n)\subseteq g(B_1)\cup \cl^{k_1}(B',
\cM_n)$
\end{enumerate}
\item We may write $B_0$ instead $\langle B_0, B_1\rangle$ if $B_1=N$
\item We say ``normally simply$^*$ good" if (iv) is replaced by
\begin{enumerate}
\item[(iv)']    $\cl^{k_1}(B',\cM_n)=g(\cl^k(B_0, N))\setminus B$. 
\end{enumerate}
\end{enumerate}
\end{definition}

\begin{definition}
\label{2.21}
The 0-1 context with closure $(\GK, \cl)$ is
{\em (normally) simply$^*$ almost nice} if:
\begin{enumerate}
\item[(A)] for every $k$, $\ell$, $\ell'$ there are $m^*=m^*(k, \ell,
\ell')$, $k^*=k^*(k, \ell, \ell')$, $t=t(k, \ell, \ell')$, $k_0= k_0(k,
\ell, \ell')$, $k_1= k_1(k, \ell, \ell')$ such that for every random
enough $\cM_n$ we have

{\em if} $\bar a\in {}^\ell|\cM_n|$ and $b\in \cM_n \setminus \cl^{k^*,
m^*}(\bar a, \cM_n)$ {\em then} there are $B\subseteq \cl^{k^*,
m^*}(\bar a, \cM_n)$ and $B^*\subseteq \cM_n$ such
that
\begin{enumerate}
\item[$(\alpha)$] $|B|\leq t$, $\bar a\subseteq B$, $\cl^{k_1}(B,
\cM_n)\subseteq \cl^{k^*, m^*}(\bar a, \cM_n)$ and
\item[$(\beta)$] $B^*\supseteq B\cup [\cl^k(\bar a b, \cM_n)\setminus
\cl^{k_1}(B, \cM_n)]$
\item[$(\gamma)$] $B<_s B^*$ (so $B^*\in \cK_\infty$) {\em or} at least
there is $B'$ such that $B<_s B'$, $(B', b, \bar c) \equiv_{\ell'} (B^*,
b, \bar c)$
\item[$(\delta)$] $\nonforkin{\cM_n\restriction B^*}{\cM_n\restriction
\cl^{k_1}(B, \cM_n)}_{B^*}^{\cM_n}$
\item[$(\varepsilon)$] letting $\bar c$ list the element of $B$ and
$$\psi(\bar x)= \bigwedge\{\varphi(\bar x): \cM_n\restriction
\cl^{k_1}(\bar c, \cM_n) \models \varphi(\bar x)\mbox{ and }
q.d.(\varphi(\bar x))\leq
k_0\}$$
we have $(\cM_n \restriction B^*, \bar c, \psi(\bar x),
\bar ab, k, k_1)$ is
(normally) simply$^*$ good or at least for some $B'$, $b'$ we have
\begin{enumerate}
\item[(i)] $(B', \bar c, \psi(\bar x), \bar ab, k, k_1)$ is
(normally) simply$^*$ good
\item[(ii)] $(B^*, b, \bar c)\equiv_{\ell'}(B', b, \bar c)$
\end{enumerate}
\end{enumerate}
\end{enumerate}
\end{definition}

\begin{remark}
We may restrict $\psi$ e.g. demand that it is in $\Pi_1$ (most natural in
the cases we have.
\end{remark}

\begin{claim}
\label{2.22}
In \ref{2.14} we can replace simply by
simply$^*$, i.e.

\noindent 1) Assume $(\GK,\cl)$ is simply$^*$ almost nice. Let $\varphi(\bar
x)$ be a f.o. formula. {\em Then} for some $\psi_\varphi(\bar x)$ we have:
\begin{enumerate}
\item[\ ] for every random enough $\cM_n$ and $\bar a\in {}^{\lg(\bar x)}
|\cM_n|$
\item[($\ast$)] $\quad\cM_n\models\varphi(\bar a)$ if and only if $\cM_n
\upharpoonright\cl^{k_\varphi}(\bar a,\cM_n)\models\psi_\varphi(\bar a)$
\end{enumerate}

\noindent 2) We have $[ \varphi \in \Pi_n \Rightarrow \psi_\varphi \in \Pi_n]$,
$[ \varphi \in \Sigma_n \Rightarrow \psi_\varphi \in \Sigma_n]$.
\end{claim}

\begin{conclusion}
\label{2.25}
\begin{enumerate}
\item Assume that the 0-1 context with closure $(\GK, \cl)$ is (normally)
simply$^*$
almost nice. Then $\GK$ satisfies the 0-1 law {\em iff} for any $k$,
$m$ we have $\langle \cM_n\restriction \cl^{k, m}(\emptyset): n<
\omega\rangle$ satisfies the 0-1 law.
\item Assume $(\GK, \cl)$ is simply$^*$ almost nice. Then $\GK$ has
convergence very weak $0-1$ law iff for every $k$, $m$ $\langle
\cM_n\restriction \cl^{k,
m}(\emptyset): n<\omega\rangle$ satisfies convergence very weak $0-1$ law.
\end{enumerate}
\end{conclusion}

\section{Further abstract closure context}
The context below is not used later so it can be skipped but it seems natural.
In this section we are lead to deal with the $0-1$ law holding for monadic
second order logic (i.e. we quantify over the sets). For this aim we will
use similar tools to those $\S2$.
Looking again at Definition \ref{2.5} or \ref{2.5B}(2), clause (A), we note
that there is an asymmetry: we try to represent $\cl^{k,m}(\bar ab,\cM_n)$ and
some $C\subseteq\cl^{k^*,m^*}(\bar a,\cM_n)$ as free amalgamation over some
$B$, small enough (with {\em a priori} bound depending on $\lg(\bar a)$
and $k$ only, there $C=\cl^k(B, \cM_n)$). Now this basis, $B$, of free
amalgamation is included in $\cl^{k^*,m^*}(\bar a,\cM_n)$ so it is without
elements from $\cl^{k,m}(\bar ab,\cM_n)\setminus
\cl^{k^*,m^*}(\bar{a},\cM_n)$. Suppose we allow this and first we deal with
the case $\cM_n$ is a graph. Hence  a member $d$ of $\cl^{k,m}(\bar
ab,\cM_n)$ may code a subset of $\cl^{k^*,m^*}(\bar a,\cM_n)$: the set
$$
\{c\in cl^{k^*,m^*}(\bar a,\cM_n):\mbox{ the pair }\{c,d\}\mbox{ is an
edge}\}.
$$
So though we are interested in f.o. formulas $\varphi(\bar x)$ speaking on
$\cM_n$, we are drawn into having $\psi_{\varphi}(\bar x)$, the
formula speaking
on $\cl^{k_\varphi, m_\varphi}(\bar x)$, being a monadic formula. Once we
allow also three place relations and more, we have to use second order logic
(still can say which quantifiers we need because the witnesses for the
elimination will come from the extensions of the $\cl^{k,m}(\bar a,\cM_n)$).
For this elimination, thinking of
an $\cM_n$, we need that any possible kind of extension of
$\cl^{k,m}(\bar a,\cM_n)$
occurs; so in the most natural cases, $|\cl^{k,m+1}(\bar a,\cM_n)|$ may be
with $2^{|\cl^{k,m}(\bar a,\cM_n)|}$ elements, so in the natural case
 which we expect to be able to understand the situation is when there
$\cl^{k,m}(\bar a,\cM_n)<\log_*(|\cM_n|)$. Still possibly
$\cl^{k,m+1}(\bar{a},\cM_n)$ is not larger than $\cl^{k,m}(\bar{a},\cM_n)$.

However there is a big difference between the monadic
(e.g. graph where the relations coded on $\cl^{k^*, m^*}(\bar a,
\cM_n)$ by members of $\cl^k(\bar ab, \cM_n)$ are monadic)
case and the more
general case. For monadic logic addition theorems like \ref{2.6A} are known,
but those are false for second order logic.

So we have good enough reason to separate the two cases. For readability
we choose here to
generalize the ``simply almost nice with $\cK=\cK_\infty$'' case only.

\begin{context}
\label{3.0} 
As in \S 2 for $(\GK,\cl)$.
\end{context}

\begin{definition}
\label{3.1new}
1) The 0-1 context with a closure operation, $(\GK, \cl)$ is {\em s.m.a.
(simply monadically almost) nice} if $\cK=\cK_\infty$, $\cl$ is transitive 
smooth local transparent (see Definitions \ref{2.2}(3),\ref{2.3}(2),(3) and 
\ref{2.5}(4),(5)) and  
\begin{enumerate}
\item[(A)] for every $k$ and $\ell$, there are
$r=r(k,\ell)$,  $k^*=k^*(k,\ell)$ and $t_1=t_1(k,\ell)$, $t_2= t_2(k,\ell)$
such that: 

for every $\cM_n$ random enough we have:

{\em if} $\bar a\in {}^\ell(\cM_n)$, $b\in \cM_n$,
$\cl^k(\bar ab,\cM_n)\nsubseteq\cl^{k^*}(\bar a,\cM_n)$

{\em then} there are $B^*$, $B^1$, $B^2$ such that:
\begin{enumerate}
\item[$(\alpha)$] $\bar a\subseteq B^1$ and $\cl^{r}(B^1,\cM_n)\subseteq
\cl^{k^*}(\bar a,\cM_n)$ and $|B^1|\leq t_1$,
\item[$(\beta)$] $B^1\subseteq B^2$, $B^2\cap\cl^r(B^1,\cM_n)=B^1$,
$|B^2|\leq t_2$, $b\in B^2$,
\item[$(\gamma)$] $B^*\supseteq [\cl^{k}(\bar ab,\cM_n)\setminus
\cl^{r}(B^1,\cM_n)]\cup B^2$ and
$\cl^k(\bar ab, \cM_n)\subseteq B^*$ (hence
$\cl^k(\bar ab,B^*)=\cl^k(\bar ab,\cM_n)$),
\item[$(\delta)$] $\nonforkin{\cM_n\restriction B^*}{\cM_n\restriction (B^2
\cup\cl^{r}(B^1,\cM_n))}_{\cM_n\restriction B^2}^{\cM_n}$ (also here
$\nonforkin{}{}_{}^{}$ is the relation of being in free amalgamation),
\item[$(\varepsilon)$] if $Q$ is a predicate from $\tau_{\cK}$ and $\cM_n
\models Q(\bar c)$, $\bar c\subseteq\cl^{r}(B^1,\cM_n)\cup B^2$ {\em then:} 
$\Rang(\bar c)\cap B^2\subseteq B^1$ or $\Rang(\bar c)\setminus B^2$ has at
most one member; if this holds we say $B^2$ is monadic over $\cl^r(B^1,
\cM_n)$ inside $\cM_n$,
\item[$(\zeta)$] $(B^*, B^1, B^2, \bar a,b,k,r)$ is m.good (see below, m
stands for monadically), so clearly $B\in {\mathcal K}_\infty$. 
\end{enumerate}

\end{enumerate}

\noindent 2) We say $(B^*,B^1,B^2,\bar a ,b,k)$ is {\em m.good} when: $B^*$, 
$B^1$, $B^2 \in \cK_\infty$ and $B^1\leq B^2 \leq B^*$,
$\bar a\subseteq B^1,b \in B^2$ and {\em for every}
random enough $\cM_n$, and $f: B^1\hra \cM_n$, and $C^1\in\cK_\infty$ such
that $\cM_n\restriction\cl^r(f(B^1),\cM_n)\subseteq C^1$, and $f^+:B^2\hra
C^1$ extending $f$ such that $C^1=f^+(B^2)\cup\cl^r(f(B^1),\cM_n)$ and
$f^+(B^2)$  is monadic over $\cl^r(f(B^1),\cM_n)$ inside $C^1$ (see above,
but not necessarily $C^1\subseteq \cM_n$) {\em there are} $g^+:C^1
\hra \cM_n$ and $g:B^*\hra\cM_n$ such that $g\restriction B^2=
( g^+ \circ f^+)\restriction B^2$ and
$$
\nonforkin{g(B^*)}{g^+(C^1)}_{g(B^2)}^{}\quad \mbox{ and }\quad \cl^k(g(\bar
ab),\cM_n)\subseteq g(B^*) \cup \cl^r(g(B^1), \cM_n).
$$

\noindent 3) Assume $\bE\subseteq\{(C,B^1,B^2):B^1\leq B^2\leq C\in\cK\}$ is
closed under isomorphism. We say $B^2$ is $\bE$-over $D$ inside $N$ if $B^2
\leq N\in\cK$, $D\leq N$ and $(N\restriction (B^2\cup D), B^2\cap D,B^2)\in
\bE$.

\noindent 4) We say $(B^*,B^1,B^2, \bar a ,b,k,r)$ is $\bE$-good when
$B^*,B^1,B^2 \in \cK_\infty$ and $B^1\leq B^2\leq B^*$,
$\bar a\subseteq B^1$, $b \in B^2$ and
for every random enough $\cM_n$ and $f:B^1 \hra C^1
\in {\mathcal K }_\infty$ extending $f$ such that
$C^1=f^+(B^2)\cup\cl^r(f(B^1),\cM_n)$ and $f^+(B^2)$ is $\bE$-over
$cl^r(f(B^1),\cM_n)$ inside $C^1$(see above but not necessarily
$C^1\subseteq\cM_n)$ there are $g^+:C^1\hra\cM_n$ and $g:B^*\hra\cM_n$
such that $g\restriction B^2=
(g^+ \circ f^+)\restriction B^2$ and
$\nonforkin{g(B^*)}{g^+(C^1)}_{g(B^2)}^{}$ and $\cl^k(g(\bar ab),\cM_n)
\subseteq g(B^*) \cup \cl^r(g(B^1),\cM_n)$.

\noindent 5) We say $\GK$ is s.$\bE$.a nice if in \ref{3.1new}(1) we replace 
clauses $(\varepsilon)$, $(\zeta)$ by
\begin{enumerate}
\item[$(\varepsilon)'$] $B^2$ is $\bE$-over $\cl^r(B^1,\cM_n)$ inside
$\cM_n$
\item[$(\zeta)'$] ($B^*,B^1,B^2,\bar a ,b,k,r)$ is $\bE$-good.
\end{enumerate}
\noindent 6) We say $\bE$ is monadic if it is as in part (3) and
$(C,B^1,B^2) \in \bE$ implies 
$$
(\bar a \in Q^C\Rightarrow \Rang(\bar a)\cap B^2\subseteq
B^1 ) \ \vee\  ( |\Rang(\bar a)\setminus B^2|\leq 1).
$$
7) We say $\bE$ as in \ref{3.1new}(3) is simply monadic
{\em if}
 it is monadic and
for any $B^1\leq B^2\in\cK$, letting
$$
\begin{array}{ll}
\Gamma_{B^2}=\Big\{\theta(y,\bar b):& \bar b\subseteq B^2\mbox{ is with no
repetition, }\theta(y, \bar x)\mbox{ is an atomic formula},\\
\ &\qquad\qquad\qquad \mbox{ each variable actually appearing}\Big\}
\end{array}
$$
we have: the  class 
$$
\begin{array}{ll}
\Big\{(D,R_{\theta(y,\bar b)},c)_{\theta(y,\bar b)\in\Gamma_{B^2},c\in B^1}:
  & D\in\cK,\\
\ &B^1\leq D,\ R_{\theta(y,\bar b)}\mbox{ is a subset of }D\setminus B^1\mbox{
and }\\
\ &\mbox{there are }C^1, f\mbox{ such that:} (C^1,B^1,B^2) \in \bE\\
\ &D\leq C^1\in\cK,\ f:B^2\hra C^1,\ f(B^2)\cap D=B^1,\\
\ &f\restriction B^1={\rm id}_{B^1},\mbox{ and for }\theta(y, \bar b)\in
\Gamma\mbox{ we have }\\
\ &R_{\theta(y, \bar b)}=\{d\in D\setminus B^1: C^1\models \theta[d,g(\bar
b)]\}\Big\}
\end{array}
$$
is definable by a monadic formula\footnote{We can restrict ourselves to the
cases $C=\cl^k(B,C)$.}.

\noindent 8) We say that $\cl$ is monadically definable for $\cK$
{\em if} for each $k$, letting $\bar{x} = \langle x_\ell:\ell<k\rangle$
some monadic formula $\Theta_k(y,x)$ we have
$y \in \cl^k(\bar x,\cM_n)\Longleftrightarrow \cM_n \restriction
 \cl^k(\bar x,\cM_n) \models \Theta_k(\bar x,y)$ holds for every
random enough $\cM_n$.

\noindent 9) We say that ${\bf E}$ is trivial if it is $\{(C,B^1,B^2 ):
\nonfork{ C}{  B^2}_{B^1},\  B^1 \le B^2 \le C\}$.  
\end{definition}

\begin{lemma}
\label{3.2new}
Assume $(\GK,\cl)$ is s.$\bE$.a.~nice and $\bE$ is simply monadic and
$\cl$ is f.o. definable or at least monadically definable (see
\ref{3.1new}(7)). {\em Then} for every f.o.~formula $\varphi(\bar x)$ there
are $k$ and a monadic formula $\psi_\varphi(\bar x)$ such that:
\begin{enumerate}
\item[$(*)_{\varphi(\bar x)}$] for every random enough $\cM_n$, for every
$\bar a\in
{}^{\lg(\bar x)}|\cM_n|$ we have
$$
\cM_n\models \varphi(\bar a)\quad\Rightleftarrow\quad\cM_n\restriction\cl^k
(\bar{a},\cM_n)\models\psi_{\varphi}(\bar a).
$$
\end{enumerate}
\end{lemma}

\begin{discussion}
\label{3.3new}
Some of the assumptions of \ref{3.2new} are open to manipulations;
others are essential.

\noindent 1) As said above, the ``monadic'' is needed in order to use
an addition theorem (see \ref{3.4below}), the price of removing it is high:
essentially above we need that after finding the copy $g(B^2)$ realizing the
required type over $\cl^k(B^2,\cM_n)$, we need to find $g(B^*)$, or a
replacement like $B'$ in the proofs in \S 2 but only the holding of some
formula $\varphi(\ldots,b,\ldots)_{b\in B^2}$ in $B^*$ is important. Now
what if the requirements on the type of $g(B^2)$ over $\cl^r 
(B^1,\cM_n)$ are not coded by some subsets of $\cl^r (B^1,\cM_n)$
but e.g. by two place relations on $\cl^r(B^1,\cM_n)$? So naturally we allow
quantification over two place relations in the formulas
$\psi_{\varphi}(\bar x)$. But then even though
\[\nonforkin{B^*}{\cl^r(B^1, \cM_n) \cup B^2}_{B^2}^{}\]
not only the small formulas satisfied by $(B^2,b)_{b \in B^1}$ are important
but also e.g. the answer to $B^*\cong\cl^r(B^1,\cM_n)$.

It natural to demand that all possibilities for the set of small
formulas in second order logic satisfied by $B^*\cup\cl^r(B^1,\cM_n)$
occur so this may includes cases where $B^*$ has to be of cardinality much
larger than $\cl^k(\bar a,\cM_n)$. So we do not formulate such lemma. Of
course some specific information may help to control the situation.
We may however consider adding (in \ref{3.1new}), the demand:
\begin{enumerate}
\item[$\boxtimes_s$] if $Y \subseteq B^* \cup \cl^r (B^1,\cM_n)$ and $Y \cap
B^*\nsubseteq B^2$, $Y \cap \cl^r(B^1,\cM_n) \nsubseteq B^1$,\\
then $Y$ is not s-connected, that is for some $Y_1$, $Y_2$, we have $Y=Y_1
\cup Y_2$, $|Y_1\cap Y_2| \leq s$ and $\nonforkin{Y_1}{Y_2}_{Y_1\cap
Y_2}^{}$ (i.e. $\nonforkin {\cM_n \restriction Y_1}{\cM_n \restriction
Y_2}_{\cM_n \restriction Y_1 \cap Y_2}^{})$.
\end{enumerate}
In this case we can allow e.g. quantification on $2$-place relations $R$ such
that $\cM_n \restriction$ Dom$(R)$ is s-connected.

\noindent 2) If $\bE$ is monadic but not simply monadic, not much is changed:
we should allow new quantifiers in $\psi_\varphi$. Let $C^1<^{\bE}_B C^2$ if
$B\leq C^1\leq C^2$ and $(C^2,B,B\cup (C^2\setminus C^1))\in \bE$. We want the
quantifier to say for $(C^1,R_{\theta(y,\bar b)},c)_{\theta(\bar y,\bar b)\in
\Gamma,c\in B}$ that it codes $C^2$ with $C^1\leq^{\bE}_B C^2$ where
$\Gamma= \Gamma_{B\cup (C^2\setminus C^1)}$, but then the logic should be
defined such that we would be able to iterate.

The situation is similar to the case that in \S 2, we have: $\cl$ is
definable or at least monadically definable.

\noindent 3) In \ref{3.2new} we essentially demand
\begin{enumerate}
\item[$(*)$] for each $t$, for random enough $\cM_n$, for every $B\subseteq
\cM_n$, $|B|\leq t$, if $\cM_n\restriction \cl^k(\bar a,\cM_n)<^{\bE}_{\bar
a}C$ then $C$ is embeddable into $\cM_n$ over $\cl^k(\bar a,\cM_n)$.
\end{enumerate}

Of course we need this just for a dense set of such $C$'s, dense in the sense
that a monadic sentence is satisfied, just like the use of $B'$ in \ref{2.5B}.
That is we may replace clause (3) of Definition \ref{3.1new}(1)(A) by
\begin{enumerate}
\item[$(\zeta)'$] there is $B'$ such that $(B', c, b)_{c\in
B_2}\equiv_{\ell}' (B^*, c, b)_{c\in B_2}$ and $(B', B^1, B^2, \bar a,
k)$ is m.~good (and $\ell'$ large enough e.g. quantifier depth of
$\psi_{\varphi_1}$ in main case).
\end{enumerate}

\noindent 4) As we have done in \ref{2.6}(2), \ref{2.14}(2), we can add that
the number of alternation of quantifiers of $\varphi$ and the number of
(possibly) alternation of monadic quantifier of $\psi_\varphi$ are equal as
long as the depth of the formulas from ``simply monadic'' is not counted
(Always we can trivially increase the q.d. so we may ask about
$\psi_{\varphi}$ with minimal number. But for a specific $\langle \cM_n:
n<\omega \rangle$ we may get better. We can though look at minimal q.d. on
all cases then it should be trivial. 

\noindent 5) Can we find a reasonable context where the situation from
\ref{3.2new} and \ref{3.3new}(1) above holds? Suppose we draw edges as here
in $\cM_{n}^0$ and redraw in the neighborhood of each edge. Let us describe
drawing fully, this for a model on $[n]$. For each $i<j$ from $[n]$ we flip
a coin $\cE_{i,j}$ on whether we have $(i,j)$ as a pre-edge, with
probability $p^n_{i,j}$. If we succeed for $\cE_{i,j}$ then for any pair
$(i',j')$ from $[n]$ we flip a coin $\cE_{i,j,i',j'}$ with probability
$p^n_{i,j,i',j'}$. The flippings are independent and finally for $i'<j'$,
$(i',j')$ is an edge if and only if for some $i<j$, $(i,j)$ is a pre-edge,
that is we succeed in $\cE_{i,j}$ and
we also succeed in $\cE_{i,j,i',j'}$. For our case let ($\alpha\in
(0,1)_{\bbr}$ is irrational):
\medskip

\noindent{\sc Distribution 1}
$$
p^n_{i,j}=p_{|i-j|}=\left\{
\begin{array}{ll}
1/|i-j|^\alpha &\mbox{when } |i-j|>1\\
1/2^\alpha            &\mbox{if\ \ } |i-j|=1\\
\end{array}\right.
$$
and $p^n_{i,j,i',j'}=\frac{1}{2^{|i-i'|+|j-j'|}}$;
\medskip

\noindent{\sc Distribution 2}

\noindent $p^n_{i,j}$ is as above and
$$
p^n_{i,j,i',j'}=\left\{
\begin{array}{ll}
\frac{1}{2^{|i-i'|+|j-j'|}}&\mbox{if } i=i'\ \vee\ j=j'\\
0                          &\mbox{if otherwise}.
\end{array}\right.
$$
\smallskip

\noindent Now distribution 2 seems to give us an example as in Lemma
\ref{3.2new}, distribution 1 fits the non-monadic case. Distribution $1$
will give us, for some pre-edges $(i,j)$, a lot of edges in the
neighbourhood of it; of course for the average pre-edge there will be
few. This give us a lot of $\leq_i$ extensions in that neighbourhood. We may
wonder whether actually the 0--1 law holds. It is intuitively clear that for
distribution 2 the answer is ``yes'', for distribution 1 the answer is
``no''. 

\noindent 6) Why in distribution 1 from (5) the 0--1 law should fail
(in fact fails badly)?  It seems to me that for distribution 1 we can find
$A\subseteq B$ such that for every random enough $\cM_n$, for some
$f:A\hra\cM_n$, the number of $g:B\hra\cM_n$ extending $f$ is quite large,
and on the set of such $g$ we can interpret an initial segment $N_f$ of
arithmetic even with $f(A)$ a segment, $N_f$ in its neighbourhood. The
problem is to compare such $N_{f_1}, N_{f_2}$ with possibly distinct
parameters, which can be done using a path of pre-edges from $f_1(A)$ to
$f_2(A)$. But this requires further thoughts.

\noindent The case of distribution $2$ should be similar to this paper.

We intend to return to this.

\noindent 7) In ${ \bE}$ is trivial, then the claim above becomes (a variant
) of the main claims in section 2 (the variant fulfill promises there).
\end{discussion}

\Proof of \ref{3.2new}:

This proof is similar to that of Lemma \ref{2.6} and \ref{2.14}.
We say in the claim that $\psi_{\varphi}(\bar x)$ or
$\psi_{\varphi}(\bar x)$, $k_{\varphi}$ witness $(\ast)_{{\varphi}(\bar x)}$.
We prove the statement by induction on q.d.$(\varphi(\bar x))$ and first
note (by clause (d) of Definition \ref{2.1}) that $(\ast)_{{\varphi}(\bar x)}
\Longrightarrow (\ast)^+_{{\varphi}(\bar x)}$ where $\psi_{\varphi}(\bar x)$
will be monadic logic.
\begin{enumerate}
\item[$(\ast)^+_{{\varphi}(\bar x)}$] for every random enough $\cM_n$,
for every $\bar a \in {}^{lg(\bar x)}(\cM_n)$ and $N$ if $\cM_n \restriction
\cl^{k_\varphi} (\bar a,\cM_n)\subseteq N\subseteq\cM_n$ then $\cM_n \models
\varphi [\bar a] \Longleftrightarrow N \models \psi_{\varphi}[\bar a]$.
\end{enumerate}

\noindent Case 1: Let $\varphi(\bar x)$ be an atomic formula. Trivial.

\noindent Case 2: $\varphi(\bar x)$ a Boolean combination of atomic formulas
and formulas $\varphi(\bar x)$ of the form $\exists y \varphi'(\bar x,y)$,
$\varphi'$ of quantifier depth$<r$, such that
$(\ast)_{\exists y \varphi'(\bar x,y)}$ holds.
Clearly follows by case 3 and case 1.

\noindent Case 3: $\varphi(\bar x)=(\exists y)
\varphi_1 (\bar x,y)$. Let $k_{\varphi_1}$, $\psi_{\varphi_1}$ be a
witness for $(\ast)_{\varphi(\bar x)}$
of \ref{3.2new} and let $k_{\varphi_1}$ $\psi^2_{\varphi_1}$ be
witness for $(\ast)^+_{{\varphi_1}(\bar x)}$ holds for it (for $\varphi_1$).
Let $r=r(k_{\varphi_1},\lg(\bar x))$,
$k^*=k^*(k_{\varphi_1},\lg(\bar x))$, $t_1=t_1(k,\lg(\bar x))$
and $t_2=t_2(k, \lg \bar x)$ be as in Definition \ref{3.1new}(1)(A), more
exactly its \ref{3.1new}(3) variant.
Let $k_{\varphi}$ be $k^*$.

It is enough to prove the following two statements:
\medskip

\noindent{\em Statement 1:}  There is $\psi^{1}_\varphi(\bar x)$ a monadic
formula such
that:
\begin{enumerate}
\item[$(\ast)_1$] for every random enough $\cM_n$, for every $\bar a\in
{}^{\lg(\bar x)} |\cM_n|$ we have $(\alpha)_1\Leftrightarrow (\beta)_1$ where:
\begin{enumerate}
\item[$(\alpha)_1$] $\cM_n\upharpoonright\cl^{k^*}(\bar a,\cM_n)\models
\psi^1_\varphi(\bar a)$
\item[$(\beta)_1$] $\cM_n\models$``there is $b$ satisfying
$\cl^{k_{\varphi_1}}(\bar{a}b,\cM_n)\subseteq\cl^{k^*}(\bar a,\cM_n)$
such that $\varphi_1(\bar a,b)$ holds.''
\end{enumerate}
\end{enumerate}
\medskip

\noindent{\em Statement 2:} There is $\psi^2_\varphi(\bar x)$ a monadic
formula such that:
\begin{enumerate}
\item[$(\ast)_2$] for every random enough $\cM_n$ and for every $\bar a\in
{}^{\lg(\bar x)}|\cM_n|$ we have $(\alpha)_2\Leftrightarrow(\beta)_2$ where:
\begin{enumerate}
\item[$(\alpha)_2$] $\cM_n\upharpoonright\cl^{k^*}(\bar a,\cM_n)\models
\psi^2_\varphi(\bar a)$
\item[$(\beta)_2$] $\cM_n\models$ ``there is $b$ satisfying
$\cl^{k_{\varphi_1}}(\bar{a}b,\cM_n)\nsubseteq\cl^{k^*}(\bar a,\cM_n)$
such that $\varphi_1(\bar a,b)$ holds''
\end{enumerate}
\end{enumerate}
(note: $(\beta)_1$, $(\beta)_2$ are complementary, but it is enough that
always at least one holds).

\noindent Note that as ``$y\in\cl^{k^*}(\bar x)$'' is monadically
definable, by \ref{3.1new}(7) clause (d) and by the choice of $k_\varphi$ we
can in $(\alpha)_2$ replace $\cl^{k^*}$ by $\cl^{k_\varphi}$, changing
$\psi^2_\varphi$ to $\psi^{2.5}_{\varphi}$, and similarly in $(\alpha)_1$
replace $\cl^{k^*}$ by  $\cl^{k^*_\varphi}$ changing $\psi^1_\varphi$ to
$\psi^{1.5}_\varphi$. 

\noindent Clearly these two statements are enough and $\psi^{1.5}_\varphi
(\bar x)\vee\psi^{2.5}_\varphi(\bar x)$ is as required.
\medskip

\noindent{\em Proof of statement 1:}

\noindent Easily, by the induction hypothesis and by the fact that the closure
is sufficiently definable.
\medskip

\noindent{\em Proof of statement 2:}

\noindent We will use a series of equivalent statements $\otimes_\ell$.
\begin{enumerate}
\item[$\otimes_1$] is $(\beta)_2$,
\item[$\otimes_2$] there are $b$, $B$ and $B^*$, $B_1$, $B_2$ such that:

$b\in\cM_n$, $\cl^{k_{\varphi_1}}(\bar{a}b,\cM_n)\nsubseteq\cl^{k^*}
(\bar a,\cM_n)$, $\bar a\subseteq B_1\subseteq\cl^{k^*}(\bar a,\cM_n)$,

$\cl^r(B_1,\cM_n)\subseteq\cl^{k^*}(\bar a,\cM_n)$, $|B_1|
\leq t_1$, $|B_2|\leq t_2$, $B_1\leq B_2\leq B^*$,

$b\in B^*$, $B^*\setminus B_1$ disjoint to $\cl^r(B_1,\cM_n)$ and
\footnote{ the $B'$ does not appear for simplicity only}
$B_1\leq_s B^*\in\cK_\infty$ and
$\nonforkin{B^*}{\cl^r (B_1,\cM_n)}_{B_2}^{\cM_n}$ and
$\cl^{k_{\varphi_1}}(\bar ab, \cM_n)\subseteq B^*$ (hence
$\cl^{k_{\varphi_1}}(\bar ab,B^*)=\cl^{k_{\varphi_1}}(\bar ab,\cM_n)$) and

$(B^*,B_1,B_2,\bar a,b,k,r)$ is $\bE$--good and
\begin{enumerate}
\item[\ $\oplus_2$] $\cM_n\models\varphi_1(\bar a,b)$.
\end{enumerate}
\end{enumerate}

\begin{enumerate}
\item[$(\ast)_2$] $\otimes_1 \Leftrightarrow \otimes_2$
\end{enumerate}
Why?  The implication $\Leftarrow$ is trivial, the implication $\Rightarrow$
holds by clause (A) in the definition \ref{3.1new}.

\begin{enumerate}
\item[$\otimes_3$] like $\otimes_2$ but replacing $\oplus_2$ by
\begin{enumerate}
\item[$\oplus_3$] $\cM_n\upharpoonright\cl^{k_{\varphi_1}}(\bar ab,M_n)\models
\psi_{\varphi_1}(\bar a,b)$.
\end{enumerate}
\end{enumerate}

\begin{enumerate}
\item[$(\ast)_3$]  $\otimes_2 \Leftrightarrow \otimes_3$
\end{enumerate}
Why?  By the induction hypothesis i.e. choice of $k_{\varphi_1}$
$\psi_{\varphi_1}$.

\begin{enumerate}
\item[$\otimes_4$] like $\otimes_3$ replacing $\oplus_3$ by
\begin{enumerate}
\item[$\oplus_4$] $\cM_n\upharpoonright[B^* \cup
\cl^{k_{\varphi_1}}(B,\cM_n)]\models\psi^2_{\varphi_1}(\bar a,b)$.
\end{enumerate}
\end{enumerate}

\begin{enumerate}
\item[$(\ast)_4$] $\otimes_3 \Leftrightarrow \otimes_4$
\end{enumerate}
Why? By $(\ast)^+_{\varphi_1}$ being witnessend by $\psi^2_{\varphi_1}$,
$k_{\varphi_1}$ see the beginning of the proof, the definition of $B^*$ and
the choice of $\psi^2_{\varphi_1}$.

For notational simplicity we assume $B\neq\emptyset$, and similarly assume
$\bar a$ is with no repetition and apply Lemma \ref{3.4below} below with the
vocabulary $\tau_{\cK}$ to the case $s=\ell$, $\bar z^2$ empty, $\bar
z^1=\langle z^1_1\rangle$, $\bar z=\langle z_1,\ldots,z_\ell\rangle$,
and $\psi(\bar z,\bar z^1,\bar z^2)=\psi(\bar
z,z^1_1)=\psi_\varphi^2(\langle z_1, \ldots, z_{\lg \bar x}\rangle,
z^1_1)$ and get $i^*$, $\theta^1_{i}(\bar z, \bar z^1)$
and $\theta^2_{i}(\bar z)$ for $i<i^*$ as there.

\noindent Next let
\begin{enumerate}
\item[$\otimes_5$] like $\otimes_4$ but replacing $\oplus_5$ by
\begin{enumerate}
\item[$\oplus_5$] letting $c_1,\ldots,c_{t_2}$ list $B_2$ possibly with
repetitions but such that $\{c_1,\ldots,c_{t_1}\}=B_1$ and $\langle c_1,
\ldots,c_{\lg(\bar x)}\rangle=\bar a$ and {\em there is} $i<i^*$ such that:
\begin{enumerate}
\item[(i)]  $B^*\models\theta^1_{i}[\langle c_1,\ldots,c_{t_2}\rangle,b]$
\item[(ii)] $\cM_n\restriction(B_2\cup\cl^k(B_1,\cM_n))\models
\theta^2_{i}[\langle c_1,\ldots,c_{t_2}\rangle]$.
\end{enumerate}
\end{enumerate}
\end{enumerate}
Now
\begin{enumerate}
\item[$(*)_5$] $\otimes_4 \Leftrightarrow \otimes_5$
\end{enumerate}
Why? by the choice of $\theta^1_{i}$, $\theta^2_{i}$ ($i<i^*$).

Let $\cP=\{(N,c_1,\ldots,c_{t_2}): N\in\cK_\infty$, with the set of elements
$\{c_1,\ldots,c_{t_2}\}\}$. Let $\{(N_j,c_1^j,\ldots,c^j_{t_2}): j<j^*\}$ list
the members of $\cP$ up to isomorphism, so with no two isomorphic. For every
$j<j^*$ and $i<i^*$ choose if possible $(N_{j,i},c^j_1,\ldots,c^j_{t_2},b^j_i)$
such that:
\begin{enumerate}
\item[(i)]   $N_j\leq_s N_{j,i}$ (in $\cK_\infty$),
\item[(ii)]  $b^j_i\in N_{j,i}\setminus N_j$,
\item[(iii)]  $N_{j,i}\models\theta^1_{i}(\langle
c^j_1,\ldots,c^j_{t_2}\rangle, b^j_i)$ and
\item[(iv)]   $(N_{j,i},\{c^j_1,\ldots,c^j_{t_1}\},\{c^j_1,\ldots,c^j_{t_2}\},
\{c^j_1,\ldots,c^j_{\lg\bar x},b^j_i\},k)$ is $\bE$--good.
\end{enumerate}
Let
$$
w=\{(i,j):i<i^*,j<j^*\mbox{ and }(N_{j,i},c^j_1,\ldots,c^j_t,b^j_i)\mbox{ is
well defined}\}.
$$
Let $\Gamma=\{\theta(y,\bar{x}):\theta\mbox{ is a basic formula,
}\bar{x}\subseteq \{x_1,\ldots,x_{t_2}\}\}$.

As $\bE$ is simply monadic (see Definition \ref{3.1new}(4)) we have: for some
monadic formula $\theta^3_i$ such that
\begin{enumerate}
\item[(*)] if $\{d_1,\ldots,d_{t_1}\}\leq C\in\cK$
letting  $ \Gamma =^{\rm df }
\{ \theta(y, \dots , x_{i(\ell) } , \dots )_{ \ell \ell(*) }  :
 \theta $ an atomic formula
for $ \tau_{\mathcal K} $,
every variable actually appear
and $ i ( \ell ) \in \{ 1, \dots , t_2\} \}  $ ;

the following are {\em equivalent}:

(a)\quad there are subsets $ R_\theta $ of $ C $  for $ \theta \in \Gamma $
and  there are $ C_1, d_t (t=t_1 +1,\dots, t_2 )$ satisfying
$R_{\varphi(y,\bar{x})}\subseteq C$ and $C\leq C_1\in\cK$, 
$C_1\setminus C=\{d_{t_1+1},\ldots, d_{t_2}\}$, and
$$
R_{\theta(y,\ldots,x_i,\ldots)}=\{e\in C:C_1\models\theta[e,\ldots,d_i,\ldots]
\}\quad\mbox{ for }\theta(y,\ldots,x_i,\ldots)\in\Gamma
$$
and $C_1\models\theta^2_i[d_1,\ldots,d_{t_2}]$

(b)\quad $C\models\theta^3_i[d_1,\ldots,d_{t_1}]$.
\end{enumerate}

Let
\begin{enumerate}
\item[$\otimes_6$] there are $b$, $B_1$ such that: $b\in\cM_n$,
$\cl^{k_{\varphi_1}}(\bar{a}b,\cM_n)\nsubseteq\cl^{k^*}(\bar a,\cM_n)$,
$\bar a\subseteq B_1\subseteq\cl^{k^*}(\bar a,\cM_n)$,
$\cl^r(B_1,\cM_n)\subseteq \cl^{k^*}(\bar a,\cM_n)$, $|B|\leq
t_1(k_{\varphi_1},\lg(\bar x))$, and
\begin{enumerate}
\item[$\oplus_6$] for some $c_1,\ldots,c_{t_1}$ listing $B_1$ such that $\bar
a=\langle c_1,\ldots,c_{\lg\bar x}\rangle$ {\em there are} $i<i^*$, $j<j^*$
such that:
\begin{enumerate}
\item[(i)]  $(\cM_n\restriction B^1,c_1,\ldots,c_{t_1})\cong (N_j,c^j_1,\ldots,
c^j_{t_1})$ i.e. the mapping $c^j_1\mapsto c_1$, $c^j_2\mapsto c_2$ embeds
$N_j$ into $\cM_n$,
\item[(ii)] $\cM_n\restriction\cl^{k_{\varphi_1}}(B^1,\cM_n)\models
\theta^3_{i}(\langle c_1,\ldots,c_{t_1}\rangle)$.
\end{enumerate}
\end{enumerate}
\end{enumerate}

\begin{enumerate}
\item[$(*)_6$] $\otimes_5 \Leftrightarrow \otimes_6$.
\end{enumerate}
Why? For proving $\otimes_5 \Rightarrow \otimes_6$ let $c_1,\ldots, c_t$ as
well as $i<i^*$ be as in $\oplus_5$, let $j<j^*$ be such that $(\cM_n
\restriction B^1, c_1,\ldots,c_{t_1})\cong (N_j,c^j_1,\ldots,c^j_{t_1})$.
A main point
is that $B^*$ exemplifies that $(i,j)\in w$.\\
For proving $\otimes_6\Rightarrow\otimes_5$ use part (B) of Definition
\ref{2.5}.

Now we have finished as $\otimes_6$ can be expressed as a monadic formula
straightforwardly. So we have carried the induction hypothesis on the
quantifier depth thus finishing the proof. \hfill\rqed$_{\ref{3.2new}}$

The following is the parallel of \ref{2.6A} for monadic logic (see Gurevich
\cite{Gu}, more \cite{Sh:463}).

\begin{lemma}
\label{3.4below}
For finite vocabulary $\tau$ and monadic formula (in
the vocabulary $\tau$) $\psi(\bar z,\bar z^1,\bar z^2)$, $\bar z=\langle
z_1,\ldots,z_s\rangle$, there are $i^*\in\bbn$ and monadic $\tau$- formulas
$\theta^1_i(\bar z, \bar z^1)=\theta^1_{i,\psi}(\bar z,\bar z^1)$,
$\theta^2_i(\bar x,\bar z)=\theta^2_{i,\psi}(\bar z,\bar z^2)$
for $i<i^*$ each of quantifier depth at most that of $\psi$ such that:
\begin{quotation}
\noindent {\em if} $\nonforkin{N_1}{N_2}_{N_0}^{N}$, $N_1\cap N_2= N_0$,
$N_1\cup N_2=N$ and the set of elements of $N_0$ is $\{c_1,\ldots,c_s\}$,
$\bar c=\langle c_1,\ldots,c_s\rangle$ and $\bar c^1\in{}^{\lg\bar z^1}(N_1)$
and $\bar c^2 \in {}^{\lg \bar z^2}(N_2)$

\noindent {\em then}
$$
\hspace{-1.2cm}N\models\psi[\bar c,\bar c^1,\bar c^2]\quad\mbox{ iff
}\quad\mbox{ for some }i<i^*, N_1\models\theta^1_i[\bar c,\bar c^1]\mbox{ and
} N_2\models\theta^2_i [\bar c,\bar c^2].
$$
\end{quotation}
\end{lemma}


\begin{thebibliography}{LuSh 435}
\makeatletter \renewcommand{\@biblabel}[1]{[#1]} \makeatother
\def\eprintfootnotetext{References of the form {\tt math.XX/$\cdots$}
 refer to the {\tt xxx.lanl.gov} archive}
\ifx\documentstyle\undefinedcontrolsequence
   \def\anyfootnote{\footnote{*}}
   \else\def\anyfootnote{\footnote}\fi
\def\eprintfn{\ifEprint\anyfootnote{\eprintfootnotetext}\fi\Eprintfalse }
\newif\ifEprint  \Eprinttrue

\bibitem[Bl96]{Bl96}John~T. Baldwin.
\newblock Near model completeness and 0--1 laws.
\newblock {\em preprint}, 1996.

\bibitem[BlSh 528]{BlSh:528}John~T. Baldwin and Saharon Shelah.
\newblock {Randomness and Semigenericity}.
\newblock {\em {Transactions of the American Mathematical Society}}, {\bf
  349}:1359--1376, 1997.
\newblock math.LO/9607226.

\bibitem[BoSp]{BoSp}Ravi~B. Boppana and Joel Spencer.
\newblock Smoothness laws for random ordered graphs.
\newblock In {\em Logic and random structures (New Brunswick, NJ, 1995)},
  volume~33 of {\em DIMACS Ser. Discrete Math. Theoret. Comput. Sci.}, pages
  15--32. American Mathematical Society, Providence, Rhode Island, 1997.

\bibitem[CK]{CK}Chen~C. Chang and Jerome~H. Keisler.
\newblock {\em {Model Theory}}, volume~73 of {\em Studies in Logic and the
  Foundation of Math.}
\newblock North Holland Publishing Co., Amsterdam, 1973.

\bibitem[Gu]{Gu}Yuri Gurevich.
\newblock {Monadic Second--Order Theories}.
\newblock In J.~Barwise and S.~Feferman, editors, {\em {Model Theoretic
  Logics}}, Perspectives in Mathematical Logic, chapter XIII, pages 479--506.
  Springer-Verlag, New York Berlin Heidelberg Tokyo, 1985.

\bibitem[LuSh 435]{LuSh:435}Tomasz {\L}uczak and Saharon Shelah.
\newblock {Convergence in homogeneous random graphs}.
\newblock {\em {Random Structures \& Algorithms}}, {\bf 6}:371--391, 1995.
\newblock math.LO/9501221.

\bibitem[Sh 550]{Sh:550}Saharon Shelah.
\newblock {0--1 laws}.
\newblock {\em {Preprint}}.
\newblock math.LO/9804154.

\bibitem[Sh 637]{Sh:637}Saharon Shelah.
\newblock {0.1 Laws: Putting together two contexts randomly }.
\newblock {\em {in preparation}}.

\bibitem[Sh 581]{Sh:581}Saharon Shelah.
\newblock {When 0--1 law hold for $G_{n,\bar{p}}$, $\bar{p}$ monotonic}.
\newblock {\em {in preparation}}.

\bibitem[Sh 517]{Sh:517}Saharon Shelah.
\newblock {Zero one laws for graphs with edge probabilities decaying with
  distance. Part II}.
\newblock {\em {Fundamenta Mathematicae}}, {\bf submitted}.

\bibitem[Sh 463]{Sh:463}Saharon Shelah.
\newblock {On the very weak $0-1$ law for random graphs with orders}.
\newblock {\em {Journal of Logic and Computation}}, {\bf 6}:137--159, 1996.
\newblock math.LO/9507221.

\bibitem[Sh 548]{Sh:548}Saharon Shelah.
\newblock {Very weak zero one law for random graphs with order and random
  binary functions}.
\newblock {\em {Random Structures \& Algorithms}}, {\bf 9}:351--358, 1996.
\newblock math.LO/9606230.

\bibitem[ShSp 304]{ShSp:304}Saharon Shelah and Joel Spencer.
\newblock {Zero-one laws for sparse random graphs}.
\newblock {\em {Journal of the American Mathematical Society}}, {\bf
  1}:97--115, 1988.

\bibitem[Sh:F192]{Sh:F192}{Shelah, Saharon}.
\newblock {Lecture notes on 0--1 laws, October'95, Rutgers}.

\bibitem[Sp]{Sp}Joel Spencer.
\newblock Survey/expository paper: zero one laws with variable probabilities.
\newblock {\em {Journal of Symbolic Logic}}, {\bf 58}:1--14, 1993.

\end{thebibliography}

\def\germ{\frak} \def\scr{\cal} \ifx\documentclass\undefinedcs
  \def\bf{\fam\bffam\tenbf}\def\rm{\fam0\tenrm}\fi 
  \def\defaultdefine#1#2{\expandafter\ifx\csname#1\endcsname\relax
  \expandafter\def\csname#1\endcsname{#2}\fi} \defaultdefine{Bbb}{\bf}
  \defaultdefine{frak}{\bf} \defaultdefine{mathfrak}{\frak}
  \defaultdefine{mathbb}{\bf} \defaultdefine{mathcal}{\cal}
  \defaultdefine{beth}{BETH}\defaultdefine{cal}{\bf} \def\bbfI{{\Bbb I}}
  \def\mbox{\hbox} \def\text{\hbox} \def\om{\omega} \def\Cal#1{{\bf #1}}
  \def\pcf{pcf} \defaultdefine{cf}{cf} \defaultdefine{reals}{{\Bbb R}}
  \defaultdefine{real}{{\Bbb R}} \def\restriction{{|}} \def\club{CLUB}
  \def\w{\omega} \def\exist{\exists} \def\se{{\germ se}} \def\bb{{\bf b}}
  \def\equivalence{\equiv} \let\lt< \let\gt> \def\implies{\Rightarrow}

\shlhetal
\end{document}